\documentclass[10pt]{amsart}
\usepackage{amsmath, amssymb, amscd, amsthm, euscript,xypic}
\newtheorem{thm}{Theorem}[section]
\newtheorem{cor}[thm]{Corollary}
\newtheorem{prop}[thm]{Proposition}
\newtheorem{lem}[thm]{Lemma}
\theoremstyle{definition}
\newtheorem{defn}[thm]{Definition}
\newtheorem{exam}[thm]{Example}
\newtheorem{rem}[thm]{Remark}

\newcounter{llistadepth}

\newenvironment{manlist}[1]{\addtocounter{llistadepth}{1}
      \edef\llistacontador{llista\romannumeral\the\value{llistadepth}}
      \list{({#1{\llistacontador}})}{\usecounter{\llistacontador}
      \def\makelabel##1{\hss\llap{##1}}
      \itemsep=2pt\parsep=0pt\topsep=3pt plus 1pt minus 1 pt}}{\endlist
      \addtocounter{llistadepth}{-1}}

\newenvironment{romlist}{\begin{manlist}{\roman}}{\end{manlist}}

\renewcommand{\H}{H^{\tiny S^1}}
\renewcommand{\cH}{H_{\tiny S^1}}
\renewcommand{\L}{\mathcal{L}}

\renewcommand{\o}{\mathcal{O}}
\newcommand{\g}{\mathfrak{g}}
\renewcommand{\u}{U\g}
\newcommand{\gl}{\mathfrak{gl}}

\newcommand{\R}{\mathbb{R}}
\newcommand{\F}{\tiny \mathbb{F}}
\renewcommand{\Re}{\mathrm{Re}\text{ }}
\newcommand{\tr}{\mathrm{tr}}
\newcommand{\Gl}{\mathrm{GL}}
\newcommand{\Sl}{\mathrm{SL}}
\newcommand{\Sp}{\mathrm{Sp}}
\newcommand{\dt}{d_{\nabla}}
\newcommand{\dvt}{d_{ev^\ast \nabla}}
\renewcommand{\O}{\mathrm{O}}
\newcommand{\C}{\mathbb{C}}

\newcommand{\MC}{\mathrm{MC}}

\newcommand{\adp}{\mathrm{ad}P}
\newcommand{\coad}{\mathrm{ad}P^\ast}
\newcommand{\adu}{P_\mathfrak u}
\newcommand{\Adp}{P_G}
\newcommand{\hal}{\alpha}

\newcommand{\G}{\mathcal{G}}
\renewcommand{\Mc}{\mathcal{MC}}
\newcommand{\W}{\mathcal{W}}
\newcommand{\w}{\overline{\mathcal{W}}}

\newcommand{\U}{\mathrm{U}}
\renewcommand{\d}{d_{\nabla}}

\renewcommand{\l}{\langle}
\renewcommand{\r}{\rangle} \numberwithin{equation}{section}

\title{String Bracket and Flat Connections}
\author{Hossein Abbaspour  \and Mahmoud Zeinalian}

\thanks{The first author was  supported by a ChateauBriand Postdoctoral Fellowship while he was visiting \'Ecole Polytechnique 2004-2005.}

\address{Centre de Math\'ematiques, \'Ecole Polytechnique, Palaiseau
91128, France. \footnote{The first author has moved to Max-Planck Institut f\"ur Mathematik in Bonn.}}
\email{abbaspou@mpim-bonn.mpg.de}

\address{Long Island University, C.W. Post College, Brookville, NY 11548, USA.}
\email{mzeinalian@liu.edu}


\begin{document} \maketitle
\begin{abstract} Let $G \to P \to M$ be a flat principal bundle over a
compact and oriented manifold $M$ of dimension $m=2d$. We construct
a map of Lie algebras $\Psi: \H_{2\ast } (L M) \to \o (\Mc)$,
where $\H_{2\ast } (LM)$ is the even dimensional part of the
equivariant homology of $LM$, the free loop space of $M$, and
$\Mc$ is the Maurer-Cartan moduli space of the graded differential
Lie algebra $\Omega^\ast (M, \adp)$, the differential forms with
values in the associated adjoint bundle of $P$. For a
$2$-dimensional manifold $M$, our Lie algebra map reduces to that
constructed by Goldman in \cite{G2}. We treat different Lie
algebra structures on $\H_{2\ast }(LM)$ depending on the choice of
the linear reductive Lie group $G$ in our discussion. This paper provides a mathematician-friendly formulation and proof of the main result of \cite{CFP} for $G=\Gl(n,\C)$ and $\Gl(n,\R)$ together with its natural generalization to other reductive Lie groups. 
\end{abstract}

\tableofcontents

\bibliographystyle{h-elsevier2}

\section{Introduction}\label{intro}
The precursor to Chas-Sullivan string bracket \cite{CS1} was
Goldman's Lie algebra structures \cite{G2} on certain vector
spaces based on the homotopy classes of closed curves on a closed
and orientable surface $S$. The simplest one of these was defined
on the vector space $ \R \hat \pi$ generated by the set $\hat \pi$
of free homotopy classes of closed and oriented curves on $S$.
There was a similar construction of a Lie algebra based on
unoriented curves. For the moment, let us talk about $ \R \hat
\pi$ for which the bracket of two equivalence classes of curves is
a signed summation of the curves obtained by breaking and
reconnecting two transversal representatives at each of their
intersection points one at a time. The relevance, and more
importantly the universality, of this algebraic object to geometry
was established by defining a map $\gamma \mapsto f_{\gamma}$,
from $\R \hat \pi$ to the Poisson algebra of smooth functions on
the symplectic space $Hom(\pi, G)/G$ of representations of $\pi$
into $G= \Gl (n, \C)$ or $\Gl(n, \R)$. Said slightly differently,
to a free homotopy class of a closed curve one assigns a function
on the moduli space of all flat connections modulo the gauge
group. The value of the function $f_{\gamma}$ at an equivalence
class $\alpha$ is the trace of the holonomy with respect to the
flat connection representing $\alpha$ along a oriented closed
curve representing the free homotopy class $ \gamma$. The Poisson
bracket of two such functions is identified as (see \cite{G2}),
\begin{equation*}
\{f_{\gamma}, f_{\lambda}\}= \sum_{p\in \gamma \# \lambda}
\varepsilon(p ; \gamma, \lambda)f_{\gamma_p\lambda_p},
\end{equation*}
where $\gamma_p\lambda_p$ denotes the product of the elements
$\gamma_p, \lambda_p \in \pi_1(S; p)$, and $\varepsilon(p ;
\gamma, \lambda)= \pm 1$ denotes the oriented intersection number
of $\gamma$ and $\lambda$ at $p$.

Goldman \cite{G2} showed that,
\begin{equation}\label{Goldman-bracket}
[  \gamma, \lambda]=\sum_{p\in \gamma \# \lambda} \varepsilon(p ;
\gamma, \lambda) \gamma_p\lambda_p,
\end{equation}
defines a Lie bracket on $\R \hat{\pi}$ the free vector space
generated by the conjugacy classes of $\pi$. In particular this
means that the map $\gamma \mapsto f_{\gamma}$ is a map of Lie
algebras.

Similarly, he showed that the Lie algebra structure on the vector
space based on the set of free homotopy classes of unoriented
curves corresponds to the case where $G$ is $\O(p, q)$, $\O(n,
\C)$, $\U(p,q)$, $\Sp(n,\R)$, or $\Sp(p,q)$. In this case the
Poisson bracket has the following formula,
\begin{equation*}
\{f_{\gamma}, f_{\lambda}\}= \frac{1}{2}\sum_{p\in \gamma \#
\lambda} \varepsilon(p; \gamma,
\lambda)(f_{\gamma_p\lambda_p}-f_{\gamma_p\lambda_p^{-1}}),
\end{equation*}
where once again $\gamma_p\lambda_p$ and $\gamma_p\lambda_p^{-1}$
denotes the product of $\gamma_p$ with $\lambda_p$ and its inverse
$\lambda_p^{-1}$ in $\pi_1(S; p)$, respectively. It was also
proved in \cite{G2} that,

$$
[ \gamma, \lambda]=\frac{1}{2}\sum_{p\in \gamma \# \lambda}
\varepsilon(p ; \gamma, \lambda)(
\gamma_p\lambda_p-\gamma_p\lambda_p^{-1}) ,$$ defines a Lie
bracket on $\R \hat{\pi}$ the free vector space generated by the
conjugacy classes of $\pi$.

Goldman's Lie bracket was generalized by Chas and Sullivan to a
bracket on the equivariant homology $\H_\ast (LM)$ of the free
loop space $LM$ of an oriented closed manifold $M$ of arbitrary
finite dimension $m$,
$$
[\cdot,\cdot ]:\H_i(LM)\otimes \H_j(LM)\to \H_{i+j+2-m}(LM).
$$
This makes $\H_\ast (LM)$ into a graded Lie algebra, after a shift
by $m-2$ in the grading. For an oriented surface $M$ of dimension
$m=2$ and $i=j=0$, the bracket on $\R\hat{\pi}=\H_0(LM)$ coincides
with that discovered by Goldman.

Inspired by \cite{G2} and \cite{CFP}, we do something similar for
the Chas-Sullivan bracket in this paper. More precisely, let $G$
denote $\Gl(n, \R)$ or $\Gl(n,\C)$, endowed with the invariant
function $f(g)=\Re \tr (g)$, and $\g$ its Lie algebra with the
nondegenerate invariant bilinear form $\l x,y\r=\Re \tr(xy)$. Let
$G \to P \to M $ be a principal bundle over a compact and oriented manifold $M$ of
dimension $m=2d$, with a fixed flat connection $\nabla$. We
construct a map of Lie algebras $\Psi$ from the equivariant
homology $\H_{2\ast } (LM)$ to the Poisson algebra of function on
the symplectic space $\Mc$, the formal completion of the space
$\{x\in \bigoplus_{k\geq 0}\Omega^{2k+1}(M,\adp) ~|~ \dt x+ 1/2
[x,x]=0\}/\G$  (see section \ref{section-MC} for definition) of
the differential graded Lie algebra $(\Omega^\ast  (M, \adp),
\dt)$ of differential forms with values in $\adp$, the associated
adjoint bundle of $P$ (see \cite{GG}). For a description of the
natural symplectic structure of $\Mc$ see Example \ref{ex1},
Proposition \ref{moment} and Theorem \ref{thm-GG}. Note that using
$\nabla$ as a point of reference, the Maurer-Cartan moduli space
contains a copy of the moduli space of flat connection on $G \to
P\to M$ (see the discussion in section \ref{2dim}). Here is one of  the
main theorems,

\setcounter{section}{7} \setcounter{thm}{3}

\begin{thm} For $G=\Gl (n,\C)$ or $\Gl (n,\R)$, the generalized holonomy
map,
$$\Psi:(\H_{2\ast }(LM),[\cdot,\cdot]) \to
(\o(\Mc),\{\cdot,\cdot\}),$$ is a map of Lie algebras.
\end{thm}
For a definition of generallized holonomy see section 6 (and \cite{CR}). For a loop $\lambda \in LM$, representing an element of $\H_0(LM)$, the value of the function $\Psi_{\lambda}$ at a flat connection is the trace of its holonomy along $\lambda$. The
machinery of the proof is robust enough to handle other reductive
subgroups of $\Gl(n, \C)$. In a manner similar to the discussion
in \cite{G2}, different subgroups correspond different Lie algebra
structures on the equivariant homology of free loop spaces. We
discuss the cases of $G =\O(p, q)$, $\O(n, \C)$, $\U(p,q)$,
$\Sp(n,\R)$, and $\Sp(p,q)$ (see section 7).

We provide a mathematician-friendly formulation and proof of the main result of \cite{CFP} for $G=\Gl(n,\C)$ and $\Gl(n,\R)$ together with its natural generalization to other reductive Lie groups.  The relevance of the subject matter to the Chern-Simons theory and its applications are explained  \cite{Schwarz}. 

Let us very briefly review each section. Section
\ref{StringBracket} is a short description of the loop product and
string bracket. In section \ref{Inv} we recall basic facts about
invariant functions on Lie groups and some of their byproducts.
Section \ref{section-MC} discusses the symplectic nature of the
set of all solutions of the Maurer-Cartan equation on a
differential graded algebra. We discuss the symplectic structure
of the moduli space of Maurer-Cartan equation via the process of
Hamiltonian reduction. Section \ref{GenHol} contains the main part
of the paper. The main concept here is that of the generalized
holonomy. Section \ref{Hamilt} concerns the construction of the
Lie algebra map for $G=\Gl(n, \C)$ or $\Gl(n, \R)$. Section
\ref{unor} deals with several other linear reductive Lie groups
which tie with the Lie algebra structure on the vector space based
on free homotopy classes of unoriented curves. Section \ref{2dim}
explains how for a $2$-dimensional manifold our construction
specializes to that described by Goldman in \cite{G2}. In Appendix
\ref{Vect}, we recall the moduli space of flat connections and its
relation to the moduli space of representations of the fundamental
group. We also review the basic facts about homology and
cohomology with local coefficients as well as a relevant version
of Poincar\'e duality. Appendix \ref{Appen-Hol} describes the
formula for the solutions of the time dependent linear system of
equations in terms of the Chen iterated integrals.

\thanks{\textbf{Acknowledgments:} {The authors would like to thank Jean Barge, David Chataur, Alberto Cattaneo, Jean Lannes, Riccardo
Longoni, and Dennis Sullivan for many helpful conversations and
comments. We are grateful to Carlo A. Rossi and James Stasheff who
read the first version of the paper and suggested many important
corrections and improvements. We are also thankful to Victor
Ginzburg and John Terrila for their useful suggestions regarding
the algebraic treatment of the moduli space. }

\setcounter{section}{1}

\section{String bracket}\label{StringBracket}

Let $M$ be an closed oriented manifold of dimension $m$. In
\cite{CS1,CS2}, Chas and Sullivan forged the term String Topology
by introducing various operations on the ordinary and equivariant
homologies of $LM=C^{\infty}(S^1,M)$ the free loop
space\footnote{In the literature $LM$ is  $C^0(S^1,M)$  which has
the same homotopy type of $C^{\infty}(S^1,M)$. }of $M$, where $S^1
=\R/\mathbb Z$.  The free loop space $LM$ is therefore a Fr\'echet
manifold benefiting from such tools as the differential forms,
principal bundles, connection, and others (see section 5 for more details). Throughout this paper the coefficient ring of the (co)homology is $\mathbb{Z}$ unless otherwise it is specified.

Chas and Sullivan's constructions includes a product on
$H_{*+m}(LM)$, called the \emph{loop product}, and its equivariant
version, the \emph{string bracket}, defined on the shifted
$S^1$-equivariant homology $\H_{\ast+m-2}(LM)=H_{\ast +m-2}(ES^1
\times_{S^1} LM)$. For the purposes of this paper we use a
description of the loop product found in \cite{CJ}. Let
$LM\times_M LM =\{(\gamma_1, \gamma_2) ~|~ \gamma_1(0)=\gamma_2(0)
\}\subset LM \times LM$ be the space of the pairs of loops with
the identical marked points and consider the following commutative
diagram,
\begin{equation*}
\begin{CD}
LM\times_M LM@>>> LM\times LM\\
@VVevV @VV{(ev,ev)}V\\
M@>\Delta>>M\times M
\end{CD}
\end{equation*}

Note that $LM\times_M LM$ is a codimension $m$ subspace of
$LM\times LM$ with a tubular neighborhood $ev^*(v_M)$ where $v_M$ is a normal bundle for the diagonal $M\hookrightarrow M\times M$.

 Consider the map, 
\begin{equation}\small{\tau: H_\ast(LM \times LM) \to H_\ast (Thom(ev^*(v_M)))\to
H_{*-m}(LM\times_M LM )},\end{equation} where the first map is Pontrjagin-Thom collapsing map to the Thom space of $ev^*(v_M)$, and the second is the Thom isomorphism for the normal bundle $ev^*(v_M)$. Also, consider the
map,
\begin{equation}\label{gamma-comp}
\varepsilon:H_\ast (LM\times_M LM)\to H_\ast (LM),\end{equation}
induced by the concatenation of the loops $\lambda_1$ and
$\lambda_2$ with the identical marked points, that is,
\begin{equation}
\lambda_1\circ \lambda_2(t)=\begin{cases}
    \lambda_1(2t)  & 0\leq t \leq 1/2 \\
    \lambda_2(2t-1)  & 1/2\leq t \leq 1
\end{cases}
\end{equation}

The loop product is defined as follows,
\begin{equation*}\bullet=\varepsilon\circ \tau:H_\ast (LM)\otimes H_\ast (LM)\simeq
H_\ast (LM\times LM)\to H_{*-m}(LM).
\end{equation*}

\begin{rem}
Here we have to modify the definition of $\lambda_1\circ
\lambda_2$ since the result may not be smooth. For that one has to
reparameterized the $\lambda_1$ and $\lambda_2$ in a neighborhood
of $0$ using a fixed smooth bijection of $[0,1]$ whose all
derivatives at $0$ and $1$ are zero. This is standard and does not
change the homological operations introduced above.
\end{rem}

Next, we recall the definition of the \emph{string bracket} on
$\H_{*+m-2}(LM)$ (see \cite{CS1}). Consider the degree one $S^1$-transfer map
$\mathfrak{m}_\ast :\H_{\ast}(LM)\to H_{\ast+1}(LM)$, and the degree zero map
$\mathfrak{e}_\ast :H_\ast(LM)\to \H_\ast(LM)$ induced by the
projection $LM\times ES^1 \to (LM \times_{S^1} ES^1)$ (see
\cite{CV, CKS}). The string bracket
$[\cdot,\cdot]:\H_i(LM)\otimes \H_j(LM)\to \H_{i+j+2-m}(LM)$ is
defined as follows,
\begin{equation}\label{equation-stringbracket}
[a,b]=(-1)^{|a|}\mathfrak{e}_\ast (\mathfrak{m}_\ast a\bullet
\mathfrak{m}_\ast b).
\end{equation}

It was proved in \cite{CS1} that $(\H_{*+m-2}(LM), [\cdot,\cdot])$
is a graded Lie algebra and that $\Delta=\mathfrak{m}_\ast \circ
\mathfrak{e}_\ast :H_{*+m}(LM)\to H_{*+m+1}(LM)$ together with the
loop product makes $H_{*+m}(LM)$ into a $BV$ algebra.  In fact
$\Delta$ is the map induced by the unit circle action.

Throughout this paper  $M$ is a manifold of dimension $m=2d$ and
therefore the shift in the degrees by the dimension does not
change the parity. This means $\H_{2\ast }(LM)$ may simply be
regarded as a non-graded Lie algebra. In this case the equation
(\ref{equation-stringbracket}) becomes,
\begin{equation}
[a,b]=\mathfrak{e}_\ast (\mathfrak{m}_\ast a\bullet
\mathfrak{m}_\ast b),
\end{equation}
for $a,b\in\H_{2\ast }(LM)$.

\section{Invariant functions and principal bundles}\label{Inv}
In this section, we briefly recall some basic facts and
definitions. Most of the material is taken from \cite{G2}. Let $G$
be a Lie group with Lie algebra $\g$. Assume that $\g$ is equipped
with a nondegenerate bilinear form $\l \cdot, \cdot\r$, which is
invariant, that is, $\l [x, y],z\r=\l x, [y, z]\r$, for all $x, y,
z \in \g$. We think of an element $x \in \g$ as a left invariant
derivation on $C^\infty(G)$ defined by,
\begin{equation} (x\cdot f)(g)=\frac{d}{dt}f(g \exp(tx))|_{t=0},
\end{equation}
for $x\in \g$. The universal enveloping algebra $\u=
\bigoplus_{k=0}^\infty \g^{\otimes k}/\l a\otimes b-b \otimes a-
[a,b]\r$ may then be regarded as the associative algebra of all
left invariant differential operators on $C^\infty(G)$. An
\emph{invariant function} $f:G\to \R$ is a $C^\infty$ function
which is invariant under conjugation. The \emph{variation
function} of an invariant function $f:G\to \R$, with respect to a
nondegenerate invariant bilinear form $\langle \cdot, \cdot
\rangle$, is the map $F: G \to \g$ defined by,
 \begin{equation}
\l F(g),x\r=(x\cdot f)(g).
 \end{equation}

Note that $F$ is a $G$-equivariant map with respect to the conjugation action on the domain and the adjoint action on the range.  In fact, one may extend $F: G \to \g$ to $\hat F: G \times \u \to
\g$ (using the identification $G \cong G \times \{ 1 \} \subset G
\times \u$) as follows.

\begin{equation}\label{eq-gene-invariant-function}
\hat F(g;[x_1\otimes \cdots \otimes x_k])=
\frac{\partial^k}{\partial t_1 \cdots \partial t_k}
F(g\exp(t_1x_1)\cdots \exp(t_{k}x_{k})) |_{(0, \cdots,0)}.
\end{equation}

We will later, in equations \ref{eq-general-sections-inv-func} and \ref{eq-general-sections-var-func}, apply the above to fibers of appropriate bundles.

A case of particular interest is when $G \subset \Gl(n, \C)$ is a
reductive subgroup, that is, a closed subgroup which is invariant
under the operation of conjugate transpose. It is evident that
$f(g)= \Re \tr \rho(g)$ is an invariant function. Moreover,
reductivity implies that the invariant bilinear form $\l x,y\r=\Re
\tr(\rho (x)\rho (y))$ is in fact nondegenerate. In fact, slightly
more generally, any covering of such a $G$ would enjoy the above
invariant function and bilinear form.
\begin{prop}
\label{prop-variation-function} Let $G \subset \Gl(n, \C)$ be a
reductive subgroup endowed with the invariant function $f(g)= \Re
\tr \rho(g)$. The variation function of $f$ with respect to the
nondegenerate invariant form $\l x,y\r=\Re \tr(\rho (x)\rho (y))$
is the composition,
\begin{equation*}
G\overset{\rho}{\longrightarrow}\Gl(n,\C)\overset{i}{\hookrightarrow}
\gl(n,\C)\overset{\small \text{pr}}{\longrightarrow}\g,
\end{equation*}
where $i$ is the inclusion of $n\times n$ invertible matrices in
all $n\times n$ matrices and $pr$ is the orthogonal projection
onto $\g$ with respect to $\l\cdot,\cdot\r$.
\end{prop}

\begin{cor}\label{cor-variantion-GL}\label{variation-GLN-envelop} For $G=\Gl(n,\R)$ or $\Gl(n,\C)$, endowed with the invariant function $f(g)=\Re \tr(g)$ and
invariant bilinear form $\l x,y \r=\Re \tr(xy)$, the variation
function  $F: \Gl(n,\F)\to \gl(n,\F)$ is the inclusion map  of
invertible matrices in all real and complex matrices,
respectively. Moreover, for $k\geq 1$,
\begin{equation*}
\begin{split}
\hat F(g;[x_1\otimes \cdots \otimes x_k])=gx_1 \cdots x_k.
\end{split}
\end{equation*}
\end{cor}

\begin{cor}\label{cor-SON}\label{variation-UN-envelop}Let $G =\O(p, q)$, $\O(n, \C)$, $\U(p,q)$, $\Sp(n,\R)$ and $\Sp(p,q)$, endowed with the invariant function
 $f(g)=\Re \tr(g)$ and invariant bilinear form $\l x,y\r=\Re \tr(xy)$. Then, the variation
function $F:G\to \g$ is given by $F(g)=\frac{1}{2}(g-g^{-1})$.
Moreover, for $k\geq 1$,
\begin{equation*}
\begin{split}
\hat F(g;[x_1\otimes \cdots \otimes x_k])= \frac{1}{2}g x_1\cdots
x_k+\frac{(-1)^{k+1}}{2}x_k\cdots x_1g^{-1}.
\end{split}
\end{equation*}
\end{cor}

The action of the universal enveloping algebra $\u$ on
$C^{\infty}(G)$ gives rise to an extension $\hat{f}: G \times \u
\to \R$ of $f:G\to \R$ (using the identification $G \cong G \times
\{ 1 \} \subset G \times \u$) defined by,
$$\hat{f}(g, r) = rf(g) \text{, for } r \in \R,$$ and,
\begin{equation}\label{eq-enevelop-function}
\begin{split}
\hat{f}(g;[x_1\otimes &\cdots \otimes x_{k+1}])=([x_1\otimes
\cdots \otimes x_{k+1}]\cdot f)(g)
\\&=\frac{\partial^{k+1}}{\partial t_1 \cdots \partial t_{k+1}}f(g
\exp(t_1x_1)\cdots \exp(t_{k+1}x_{k+1}))|_{(0, \cdots, 0)}.
\end{split}
\end{equation}

We have,
\begin{equation}\label{eq-general-variation}
\begin{split}
\hat{f}(g;[x_1\otimes \cdots \otimes x_{k+1}])=
\frac{\partial^k}{\partial t_1 \cdots \partial t_k}\langle \hat
F(g\exp(t_1x_1)\cdots \exp(t_{k}x_{k})),x_{k+1} \rangle|_{(0,
\cdots, 0)}.
\end{split}
\end{equation}

For the invariant function $\tr:\Gl(n, \C) \to \R$, $f(g)=\Re
\tr(g)$, and the bilinear form $\l x, y\r= \Re \tr(xy)$, we have,
\begin{equation}\label{calcul-trace}
\hat{\tr}(g;[x_1\otimes \cdots \otimes x_{k+1}])=\tr (gx_1\cdots
x_{k+1}).
\end{equation}

Let $M$ be a closed oriented manifold of dimension $m=2d$ and $G
\to P\to M$ be a principal bundle equipped with a flat connection
$\nabla$. Let $P(b)$ denote the fiber of $P$ at a point $b \in M$.
Let $\text{conj}: G \to Aut(G)$ denote the action of $G$ on itself
defined by $conj (g)( h)=g^{-1}hg$. Since conjugation fixes the
identity, $conj$ induces an action $Ad: G \to Aut(\g)$ and
subsequently an action $Ad: G \to Aut(\u)$. Let $\Adp$, $\adp$,
and $\adu$ respectively denote the associated bundles to these
representation, with fibers $\Adp (b)$, $\adp (b)$, and $\adu (b)$
at a point $b$. Then, $\Adp (b)$ is canonically identified as the
group of $G$-equivariant diffeomorphisms of the fiber $P(b)$ of
$P$ at that point $b$. This group is isomorphic to $G$, while
there is no natural choice of an isomorphism. The set of all
sections of $\Adp$, denoted by $\Gamma(\Adp)$, is then
identifiable as the group of fiber preserving and $G$-equivariant
diffeomorphisms of $P$. By the same token the adjoint
representation $Ad: G \to Aut(\mathfrak g)$ gives rise to an
associated vector bundle $\adp$ whose fiber $\adp (b)$ is a Lie
algebra canonically identified with the Lie algebra of $\Adp(b)$.
Clearly, there is a well-defined exponential map $exp: \adp \to
\Adp$ inducing and exponential map $exp: \Gamma(\adp) \to
\Gamma(\Adp)$ and there are natural actions $conj: \Adp(b) \times
\adp (b) \to \adp (b)$ and $conj: \Adp(b) \times \adu (b) \to \adu
(b)$ inducing an actions $Ad:\Gamma(\Adp) \times \Gamma(\adp) \to
\Gamma(\adp)$ and $Ad:\Gamma(\Adp) \times \Gamma(\adu)\to \Gamma
(\adu)$. Therefore there is a natural associative multiplication,
$*$, on $\Adp(b) \times \adu(b)$ for every $b$ defined as follows,
\begin{equation*}
(g, u)*(h, v)=(gh, Ad_h(u)\otimes v).
\end{equation*}
Here $Ad$ is the induced action of the group $G$ on the universal
enveloping algebra, $\u$, of its Lie algebra $\g$.

This multiplication extends naturally to $\Gamma(\Adp) \times
\Omega^\ast (M,\adu)$ that if for $(h_i,\alpha_i)\in \Gamma(\Adp)
\times \Omega^\ast(M,\adu)$, $i\in \{ 1,2 \} $,
\begin{equation}\label{star-operator}
(h_1,\alpha_1)*(h_2,\alpha_2)=\left(h_1h_2,
(Ad_{h_2}\alpha_1)\wedge \alpha_2\right),
\end{equation}
where $\wedge$ product on $\Omega(M, \adu)$ is the tensor product
of two multiplications: the exterior product of differential forms
and the multiplication of values in $\adu$. Note that unlike the
ordinary wedge product which is associative and graded
commutative, $\wedge$ is only associative.

An invariant function $f:G\to \R$ induces a function $f:\Adp\to \R
$, and subsequently $f:\Gamma(\Adp)\to C^{\infty}(M)$. Also, note
that $\hat{f}$ extends to,
 \begin{equation}\label{eq-general-sections-inv-func}
 \hat{f}:\Gamma(\Adp)\otimes \Gamma(\adu)\to C^{\infty}(M).
 \end{equation}

Similarly, $\hat F: G \times \u \to \g$ naturally extends to
 \begin{equation}\label{eq-general-sections-var-func}
 \hat{F}:\Gamma(\Adp)\otimes \Gamma(\adu)\to \Gamma(\adp).
 \end{equation}

\begin{rem}\label{remark-important}For linear groups $G=\Gl(n, \C)$ or $\Gl(n,\R)$, the Lie algebra $\g$ is actually an associative
algebra, and therefore there is a natural map of associative
algebras $\Pi:\u \to \g$, sending $\Pi:[x_1 \otimes \cdots \otimes
x_k] \mapsto x_1 \cdots x_k$. Moreover, for $\tr:G\to \R$ and $\l
x, y\r= \Re \tr(xy)$, one can show that,
\begin{equation*}
\begin{split}
\hat{\tr}((h_1,[x_1 \otimes \cdots \otimes x_k])* (h_2,[y_1
\otimes
\cdots \otimes y_l])) &=\tr(h_1\Pi[x_1 \otimes \cdots \otimes x_k]h_2\Pi[y_1 \otimes \cdots \otimes y_l]) \\
&= \tr(h_1x_1\cdots x_kh_2y_1\cdots y_l),
\end{split}
\end{equation*} where in the last line $\tr: \g \to \R$ is really the derivative of $\tr: G \to \R$ at the identity.
This gives rise to the following identity for elements $\alpha_1,
\alpha_2 \in \Omega^\ast (M, \adu)$.

\begin{equation}\label{trace}
\begin{split}\hat{\tr}((h_1,\alpha_1)*(h_2,\alpha))=\tr(h_1 \Pi\alpha_1\wedge
h_2 \Pi\alpha_2),
\end{split}
\end{equation}
where in the right hand side all multiplication are matrix
multiplications and $\tr:\Omega^\ast (M, \adp) \to C^\infty (M)$.
\end{rem}

\section{Symplectic nature of the Maurer-Cartan moduli space }\label{section-MC}

In this section, we mostly follow \cite{GG}. For the sake of
completeness and clarity we mention some of their result without
repeating their proofs.

\begin{defn}\label{def-DGLA}
A \emph{cyclic differential graded Lie algebra},
$(\L,d,[\cdot,\cdot],\omega(\cdot,\cdot))$, or a
\emph{\text{cyclic} DGLA} for short, consists of a
$\mathbb{Z}_2$-graded Lie algebra $\L= \L_0 \oplus \L_1$, with a
differential $d$, and a bilinear form $\omega:\L\times \L \to
\mathbb R$, satisfying,
\begin{romlist}
\item $d(\L_{0}) \subseteq \L_{1}$ and $d(\L_1) \subseteq \L_0$

\item $d[x, y]=[dx,y]+(-1)^{|x|}[x,dy]$ \item $\omega ([x,y],z ) =
\omega(x,[y,z])$ \item $\omega(dx,y)+(-1)^{|x|}\omega(x, d y)=0$
\item $\omega(y, x)=(-1)^{|x||y|}\omega(x, y)$ \item
$\omega(\cdot,\cdot)|_{\L_0}$ and $\omega(\cdot,\cdot)|_{\L_1}$
are nondegenerate

\item  $\omega(x, y) =0$ for $x \in \L_0$ and $y \in \L_1$
\end{romlist}

We furthermore assume that $\L_0$ is the Lie algebra of a complex,
connected, and simply-connected linear algebraic group $\G$.
\end{defn}

\begin{exam}\label{ex1}
Let $G \to P\to  M$ be a principal bundle over a compact and oriented
manifold $M$ without boundary of dimension $m=2d$ which is endowed with a flat
connection $\nabla$. Assume that the Lie algebra, $\g$, of $G$ is
equipped with a nondegenerate invariant bilinear form $\l \cdot,
\cdot\r$. Recall that invariance means $\l [x, y],z\r=\l x, [y,
z]\r$.  Consider the graded Lie algebra $\L=\Omega^\ast (M,\adp)$
of differential forms with values in $\adp$, the adjoint bundle of
$P$. In order to see the graded Lie algebra structure better,
recall that the tensor product of a graded commutative algebra and
a graded Lie algebra is a a graded  Lie algebra. Note that
$\L=\Omega^\ast (M,\adp)=\Omega^\ast (M) \otimes_{C^\infty(M)}
\Gamma(\adp)$, Here, the graded Lie algebra $\Gamma(\adp)$ is
concentrated in degree zero. The covariant derivative associated
to the connection $\nabla$ induces a derivation
$d_\nabla:\Omega^\ast (M,\adp)\to \Omega^{* +1}(M,\adp)$. Flatness
means that $\d$ is a differential, that is to say $\d^2=0$. Define
a bilinear form $\omega:\L\times\L\to \mathbb R$ as the following
composition.\begin{equation*}\small{ \omega:\Omega^i(M,\adp)\times
\Omega^j(M,\adp)\overset{\wedge}{\to} \Omega^{i+j}(M,\adp\otimes
\adp)\overset{\l\cdot,\cdot\r}{\to} \Omega^{i+j}(M,\mathbb
R)\overset{\int_M}{\to} \mathbb{R}},
\end{equation*}
if $i+j=2d$, and zero otherwise.

Clearly, condition (i) and (iii) are satisfied and condition (ii)
and (iv) follow from Corollaries \ref{App-compatible-bracket} and
\ref{app-coro-horiztonal}. It follows from the invariance property
and nondegeneracy of $\l\cdot,\cdot \r$ that conditions (v)-(vi)
hold. Lastly, (vii) is true because the manifold is even
dimensional. This means $\Omega^\ast (M, \adp)$ is a cyclic DGLA.
\end{exam}

Given a cyclic DGLA $(\L,d,[\cdot,\cdot],\omega(\cdot,\cdot))$,
clearly $(\L_1,\omega|_{\L_1})$ is a symplectic vector space. To
every $a\in \L_0$, one associates a vector field on $\L_1$ by
defining $\xi_a(x)=[a,x]-da$, for all $x \in \L_1$. The vector
field $\xi_a$ respects the symplectic structure, that is to say,
$L_{\xi_a}\omega=0$. In fact, this  infinitesimal action is
Hamiltonian and the following proposition describes its moment
map. For a discussion and proof see section 1.5 of \cite{GG}.

\begin{prop}\label{moment} Let  $(\L,d,[\cdot,\cdot],\omega)$ be a cyclic DGLA. The map $\phi:\L_1\to \L_0^\ast $ defined by
$\phi(x)(a)=\omega (dx+1/2[x,x], a)$, for all $x\in\L_1$ and
$a\in\L_0$, is a moment map for the action of $\L_0$ on the
symplectic space $(\L_1,\omega|_{\mathcal L_1})$.
\end{prop}

%\begin{proof} For $a\in\L_0$, consider $\lambda_a\in C^{\infty}(\L_1)$
%defined by $ \lambda_a(x)=\phi(x)(a)$. We must prove that in
%$\Omega^1(\L_1)$,  $d\lambda_a=i_{\xi_a} \omega$. The derivative
%of the map $x \mapsto dx + 1/2[x, x]$ at a point $x$ is given by
%$y\mapsto dy + [x, y]$, for all $y\in \L_1=T_x \L_1$, and
%therefore,
%\begin{eqnarray*}
%d_x\lambda_a(y)&=& \omega( dy+[x,y],a)\\&=&
%\omega(y,da)-\omega(y,[x,a])\\ &=& \omega(y,da-[x,a])\\&=& \omega
%(-da+[x,a],y) \\ &=&(i_{\xi_a(x)}\omega)(y)
%\end{eqnarray*}
%\end{proof}

Define,

\begin{equation*} \MC(\L)=\{x\in \L_1 ~|~ dx+1/2[x,x]=0\}.
\end{equation*}

Note that he map $ a\mapsto \xi_a$ is a Lie algebra homomorphism
and that for every $a\in \L_0$, the vector field $\xi_a$ is
tangent to $\MC(\L)$ (see Lemma 1.2.1 of \cite{GG}). This implies
a well-defined action of $\G$ on $\L_1$ by affine linear
transformations referred to as the gauge group. Note that this
action is different than the ordinary adjoint action of $\G$ on
$\L_1$ and it keeps invariant the scheme $\MC(\L)$.

%The quotient,
%\begin{equation} \label{mc}\Mc(\L)=\MC(\L)/\G
%\end{equation}
%is called the Maurer-Cartan moduli space.

Let us us be more precise about the nature of the quotient
$\MC(\L)/\G$. The space $\MC(\L)$ is a defined as the zero level
set of the moment map. Since the group $\G$ keeps $\MC(\L)$
invariant, one can consider the quotient $\MC(\L)/\G$. Note that
because $0$ is not a regular value and $\G$ does not act freely,
$\MC(\L)/\G$ is not a smooth variety. The good news is that just
as in \cite{GG} one may still treat $\MC(\L)/\G$ completely
algebraically as a formal scheme (see Remark 2.2.1 of \cite{GG}).
Thus, one starts with a formal completion $\hat\L_1$ of $\L_1$ at
the origin and considers the closed subscheme $\phi^{-1}(0)
\subset \hat\L_1$. Then, the completion of $\mathcal G$ gives a
pro-algebraic groupoid which acts on $\phi^{-1}(0)$. Passing to
the quotient yields a pro-algebraic stack $\Mc(\L)$. We use
$\o(\Mc(\L))$ to denote the coordinate ring of $\Mc(\L)$. Thus,
elements of $\o(\Mc(\L))$ have representatives in
$\o(\hat\L_1)=\mathbb C[[\hat\L_1]]$, the formal power series on
$\L_1$.  In the best of all cases one may be able to show that
certain desired series converge. This is in fact what happens in
our definition of the generalized holonomy (see (\ref{formal}) and
Definition \ref{formal2}) as explained in Remark \ref{conv}. Let
us also talk about the symplectic nature of the moduli space. In
the smooth case, where one deals with the pre-image of a regular
value of the moment map and has a free group action, the
symplectic quotient is a smooth variety with a well-defined
Zarisky tangent space at each point. Unfortunately, this is hardly
the case in almost every interesting example. In such cases, one
can talk about the tangent space $T_p \Mc$ as a $3$-term complex
$Lie(\mathcal G) \to T_p \hat\L_1 \to Lie(\mathcal G)^\ast$,
concentrated in degrees $-1$, $0$, and $1$, where $Lie(\mathcal
G)$ is the Lie algebra of $\mathcal G$ (see Section 1.6 of
\cite{GG}). It is clear that only in the nicest of all cases, when
$0$ is a regular value of the moment map and $\G$ acts freely,
that the homology of this complex is concentrated in degree $0$.
In general, one needs to think about this complex as the tangent
space (complex) to the symplectic reduction $\Mc(\L)$ at a point
$p$, keeping track of all the failures. In this enlightened view,
a symplectic structure on $\Mc(\L)$, for instance, is a
isomorphism $\omega_p: T_p\L_1 \to T_p^\ast\L_1$, establishing an
isomorphism between the tangent complex and its dual. The
symplectic structure on $\Mc$ gives $\o(\Mc(\L))$ a Poisson
algebra structure.

The following is a direct application of a theorem of Gan and
Ginzburg \cite{GG} for a reductive subgroup $G \subseteq
\Gl(n,\C)$ equipped with the nondegenerate invariant bilinear form
$\l x,y\r=\Re \tr(xy)$.

\begin{thm}\label{thm-GG} Let $G\subseteq \Gl(n, \C)$ be a reductive subgroup and $G\to P\to M$ be a principal bundle endowed with a flat connection $\nabla$.
The $2$-form $\omega(x,y)=\int _M\Re \tr(x\wedge y)$ defines a
symplectic structure on the Maurer-Cartan moduli space,
$\Mc(\Omega^\ast (M,\adp))$.
\end{thm}

For brevity, we reserve the symbol $\Mc$ for the
 Maurer-Cartan moduli space $\Mc(\Omega^\ast (M,\adp))$. Note that a typical element in $\Mc$ is a sum of differential forms of odd degrees. One
may think of $\Mc$ as an extension of the moduli space of flat
connections on $P$. In fact, if the Hard Lefschetz theorem holds
for $M$, then the moduli space of flat connections on $P$ may be
viewed as a symplectic substack of $\Mc$ (see section 1.7 of
\cite{GG} and section \ref{2dim} of this paper).

\section{Differential forms on free loop spaces}\label{diffLoop}

In this section we  shall present a model of the free loop space
$LM$ which is conducive to the notions of the de Rham differential
forms and its $S^1$-equivariant model. These differential forms
contain the Chen iterated integrals \cite{Ch1,Ch}. Moreover, the
cohomology of these forms, called the de Rham cohomology and the
equivariant de Rham cohomology respectively, compute the usual and
$S^1$-equivariant singular cohomologies of $LM$. Such a model of
the free loop space must support the usual differential geometric
notions such as connections and all related pullback diagrams.
Since the different parts of the required framework are developed
in different references, we briefly summarize some of the needed
definitions and statements and refer the reader to
\cite{H,God,Br,GS} for a fuller discussion.

%We first recall the notion of \emph{differentiable space} (possibly infinite-dimensional)
%modeled on a locally convex Hausdorff topological vector space $E$. There is a detailed  account of this  in \cite{H, Br}.
Let $E$, $F$ denote locally convex Hausdorff topological vector
space and $U\subset E$ an open set, we say that the map $f:U\to F$
is of class $C^1$ if the limit $df(x,v)=\underset{t\to
0}{\lim}\frac{f(x+tv)-f(x)}{t} $ exists and is continuous  as a
function of $(x, v)\in U\times E$. Similarly, we can define
functions of class $C^k$ and thus $C^{\infty}$. Then $\Omega^n(U)$
the space of differential $n$-forms on $U$ is defined to be the
space of smooth functions $\omega:U\times E^n\to \R$ which are
multilinear and antisymmetric in the last $n$ variables. Then, one
defines the exterior differential $d$ which satisfies $d^2=0$ (see
\cite{H, Br}). Also, the Poincar\'e lemma holds for convex open
subsets of $E$ (see \cite{Br} for a proof).

A \emph{differentiable space} modeled on $E$ is  a Hausdorff space
$N$ with a covering $\{U_i\}_{i\in I}$ and a collection of
homeomorphisms $\phi_i:U_i\to V_i\subset E$, such that the
transition maps $\phi_j\circ \phi_i^{-1}$ are smooth.  Then  a
differential $n$-form on an open $U\subset N$ is defined to be  a
collection of $\omega_i\in \Omega^n(\phi_i(U\cap U_i))$ patched
together  by the transition map $\phi_j\circ \phi_i^{-1}$. This
enables us to define the \emph{de Rham complex} of $N$ denoted $(
\Omega^\ast (N),d)$ and the \emph{de Rham cohomology} $H^\ast
_{DR}(N)$ of $N$ to be the cohomology of  $( \Omega^\ast (N),d)$.

It is known that the de Rham theorem holds (see for example
\cite{Br}). This is shown using sheaf cohomology. The de Rham
forms on the opens $U\subset N$ define a sheaf
$\underline{\Omega^\ast(N)}$ on $N$ which is a \emph{resolution}
of $\underline{\R}_N$ the sheaf of the constant functions  (by
Poincar\'e lemma). When $N$ is a paracompact differentiable space
and the sheaf $ \underline{\Omega^\ast(N)}$ is an acyclic
\emph{soft sheaf},  we can use the natural resolution $
\underline{\R}_N \to \underline{\Omega(N)}^\ast $ to calculate the
sheaf  cohomology $H_{sheaf}^\ast (N, \underline{\R}_N)$ implying
that $H^\ast _{DR}(N) \simeq H_{sheaf}^\ast (N, \underline{\R}_N)$
(see \cite{Br} 1.4.7 or \cite{God}).

As it turns out, $LM=C^{\infty}(S^1,M)$ can be made into a
\emph{paracompact} differentiable space modeled on
$C^{\infty}(S^1,\R^n)$ (p. 110 \cite{Br}). The topology of
$C^{\infty}(S^1,\R^n)$ is defined using a family of the norms $\|
\cdot\|_k$ given by

$$
\| f\|_k^2=\int_0^1(\| f(t)\|^2+\cdots+\| \frac{d^k}{dt^k}f(t)\|^2) dt,
$$
which make $C^{\infty}(S^1,\R^n)$ into a  \emph{Fr\'echet vector
space} and  $LM$ into a \emph{Fr\'echet manifold}. This allows us
to consider vector bundles with connection on $LM$ and do  the
differential geometry needed in this paper just as in the finite
dimensional case (see \cite{H} for the details).

Moreover, $C^{\infty}(S^1,\R^n)$ can be though of as an inverse
limit of the Sobolev spaces $H^{2, k}(S^1,\R^n)$, which are the
completions of $C^{\infty}(S^1,\R^n)$ with respect to the Sobolev
norms $\| \cdot \|_k$. Therefore, $C^{\infty}(S^1,\R^n)$ can be
treated as an \emph{inverse limit of Hilbert spaces}, ILH for
short.

 Because of the fact that the free loop space $LM$ is modeled on $C^{\infty}(S^1,\R^n)$ (which is an ILH) the sheaf $\underline{\Omega(N)}^\ast $ is a soft
 sheaf, and therefore,

 $$
 H^\ast _{DR}(LM)\simeq  H_{sheaf}^\ast (LM, \underline{\R}_{LM}).
 $$

It is proved in \cite{God} that for a paracompact topological
space  $N$, $$H_{sheaf}^\ast (N, \underline{\R}_N)\simeq H^\ast
(N,\R),$$ therefore,
\begin{equation}\label{eq-equi-4}
 H^\ast _{DR}(LM)\simeq H^\ast (LM ,\R).
 \end{equation}

 We turn to the equivariant de Rham theorem for which one can find in \cite{GS} a rigorous treatment for finite dimension manifolds.
 The same treatment works  for $LM$ as we have the right notion of differentials forms for $LM$ with a de Rham theorem.
 We give the steps whose proofs are exactly as same as the ones in chapter 2 of
 \cite{GS}.

Recall that for a group $G$ acting on a  space $N$,  the $G$
equivariant cohomology of N is defined to be $H^\ast
_{G}(N,\R)=H^\ast (N\times EG/G,\R)$ where $ EG$ is a contractible
topological space  on which $G$ acts freely, and  $BG= EG/G$ is
classifying of the $G$.  The case of interest for us is the action
of $S^1$ on the loop space $LM$ by changing the parametrization.
For $S^k=\{( z_1,\cdots,z_k) \in \C^k| |z_1|^2+\cdots +
|z_k|^2=1\}$ let $ES^1=S^{\infty}= \underrightarrow{\lim}S^{2k-1}$
defined by the natural inclusions $S^{2k-1} \to S^{2k+1},
(z_1,\cdots,z_k) \mapsto  (z_1,\cdots,z_k,0)$. This comes with the
natural inclusions $j_k:S^{2k-1}\to S^{\infty}$  whose opens are
the subset $U\subset S^{\infty}$ with $j_k^{-1}(U)$ is open in
$S^{2k-1}$ for all $k$.  Then $S^1$ acts freely  on $S^{\infty}$
and we have $BS^1=\C P^{\infty}=\underrightarrow{\lim}\C P^{k-1}$
whose topology is defined  similarly.  Then the real coefficient
equivariant  cohomology of  $LM$ is defined to be
 \begin{equation}\label{eq-equiv-3}
 H_{S^1}^\ast (LM,\R)= H^\ast (LM\times S ^\infty/S^1, \R),
 \end{equation}
where $S^1$ acts diagonally on $LM \times S^\infty$.

 The inclusions,
 $$
 \cdots  \to S^{2k-1} \to S^{2k+1} \to  S^{2k+3}  \to \cdots
 ,$$
 induces the sequence of the projections,
  $$
 \cdots \leftarrow \Omega^\ast (S^{2k-1}) \leftarrow  \Omega^\ast (S^{2k+1}  )\leftarrow  \Omega^\ast ( S^{2k+3} )  \leftarrow\cdots
 ,$$
which allows us to consider to the inverse limit $\Omega^\ast
(S^\infty)=\underleftarrow{\lim} \Omega^\ast (S^{2k-1})$ and we
can form the complex $(\Omega^\ast (S^\infty),d)$ which is acyclic
and satisfies condition (C) (see \cite{GS} p.
29)\footnote{Condition (C) is an algebraic way of interpreting
that an action is locally free.}. Since $S^1$ acts on $LM$ as well
as on all $ S^{2k-1}$'s,  we can consider the subcomplex of
\emph{basic forms} in $\Omega^\ast (LM)\otimes \Omega^\ast
(S^\infty)$, that is the set of all $\omega$ such that,
  \begin{equation*}
   i_X \omega=0  \text { and }  L_X\omega=0,
 \end{equation*}
where $L_X=i_X\circ d+d\circ i_X$ is the Lie derivative. We define
the equivariant de Rham cohomology to be
  $$
 H^\ast _{S^1, DR}(LM,\R)= H^\ast ((\Omega^\ast (LM)\otimes \Omega^\ast (S^\infty))_{bas}).
 $$
The rest of this section is devoted to the proof of
\begin{equation}\label{eq-equi-2}
H^\ast _{S^1,DR}(LM,\R)\simeq  H_{S^1}^\ast (LM,\R).
\end{equation}
 Consider the  diagonal action of $S^1$ on $LM\times  S^{2k-1} $ and the
 projection,
 $$
 \pi:LM\times S^{2k-1}\to LM\times S^{2k-1}/S^1.
 $$
  Since the action is free,  the subalgebra  $\pi^\ast (\Omega^\ast  (LM\times S^{2k-1} /S^1)) \subset \Omega^\ast  (LM\times S^{2k-1})$ can be
  characterized as  the basic differential forms  denoted $\Omega^\ast  (LM\times S^{2k-1})_{bas}$; moreover $\pi^\ast $ is injective.
  We denote the subcomplex of  basic forms by $\Omega^\ast (LM\times  S^{2k-1})_{bas}$ and we
  have,
  \begin{equation}\label{eq-equi-1}
   H^\ast _{DR}(LM\times  S^{2k-1}/S^1)\simeq H^\ast ( \Omega^\ast (LM\times  S^{2k-1})_{bas}).
  \end{equation}
Using the projections,
 $$
\cdots\leftarrow  \Omega^\ast (LM\times  S^{2k-1})_{bas}
\leftarrow \Omega^\ast (LM\times S^{2k+1})_{bas} \leftarrow \cdots
 ,$$
 we can form the complex  $\Omega^\ast (LM\times S^\infty)_{bas}= \underleftarrow{\lim} \Omega^\ast (LM\times  S^{2k-1})_{bas}$.

\begin{prop}
$ H_{S^1}^\ast (LM,\R)\simeq H^\ast (\Omega^\ast (LM\times
S^\infty)_{bas})$.
\end{prop}
\begin{proof}
Similarly to the proof of (\ref{eq-equi-4}), for the
differentiable space $LM\times  S^{2k-1}/S^1$, we have the
isomorphism,
$$ H^\ast _{DR}(LM\times  S^{2k-1}/S^1)=H^\ast (LM\times  S^{2k-1}/S^1,\R),$$
which is compatible the inclusions $LM\times  S^{2k-1}/S^1
\hookrightarrow LM\times  S^{2k+1}/S^1$. Therefore,
\begin{eqnarray*}
H_{S^1}^\ast (LM,\R) &\simeq&  \underleftarrow{\lim} H^\ast
(LM\times S^{2k+1}/S^1) \\ &\simeq& \underleftarrow{\lim}H^\ast
_{DR}(LM\times  S^{2k-1}/S^1)\\&\simeq& H^\ast (\Omega^\ast
(LM\times S^\infty)_{bas}).
  \end{eqnarray*}
 \end{proof}
 So, to prove (\ref{eq-equi-2}) it suffices to show that,
 \begin{equation}
 H^\ast (\Omega^\ast (LM\times S^\infty)_{bas})\simeq H^\ast _{DR}((\Omega^\ast (LM)\otimes \Omega^\ast (S^\infty))_{bas}).
 \end{equation}
 This follows from a standard spectral sequence argument
 and the fact that inclusion $\Omega^\ast (LM\times S^\infty)\hookrightarrow \Omega^\ast (LM)\otimes \Omega^\ast (S^\infty)$ induces an isomorphism in cohomology, and that
  $\Omega^\ast (S^\infty)$ is acyclic and satisfies condition (C)(see \cite{GS} p. 30).

\section{Generalized holonomy}\label{GenHol}
Let $M$ be a compact and oriented manifold without boundary. Given a principal bundle $G\to P\to M$ endowed
with a flat connection $\nabla$, the trace of the holonomy yields
a well-defined map on the set of the free homotopy classes of
loops in the base manifold $M$.  It was discussed in \cite{CR} how this may generalize to families of loops, more
precisely, to the homology classes in the free loop space of the
underlying manifold. We will give a mathematical account of this in this section.  Throughout this section we will be using the notion of differential forms  developed in section 5 for $LM$ equipped with it natural Fr\'echet manifold structure.

For every $n\geq 0$, consider the $n$-simplex,
$$\Delta^n=\{(t_0, t_1,\cdots,t_n, t_{n+1}) ~|~ 0=t_0 \leq t_1\leq \cdots\leq
t_n \leq t_{n+1}=1 \}.$$ Define the evaluation maps $ev$, and
$ev_{n,i}$, for $1\leq i\leq n$ as follows,
\begin{equation*}
\begin{split}
&ev:\Delta^n \times LM \to  M,\\
&ev(t_0,t_1, \cdots,t_n,t_{n+1}; \gamma)=\gamma(0)=\gamma (1),\\
&ev_{n,i}:\Delta^n \times LM \to  M,\\
&ev_{n,i}(t_0,t_1, \cdots,t_n,t_{n+1}; \gamma)=\gamma (t_i).
\end{split}
\end{equation*}
Let $T_i:ev_{n,i}^\ast  (\adp) \to ev^\ast  (\adp)$ denote the
map, between pullbacks of the adjoint bundles over $\Delta^n
\times LM$, defined at a point $(0=t_0, t_1, \cdots, t_n,
t_{n+1}=1; \gamma)$ by the parallel transport along and in the
direction of $\gamma$ from $\gamma(t_i)$ to
$\gamma(t_{n+1})=\gamma(1)$, in $\adp$ with respect to the flat
connection $\nabla$. Note that if $R$, $S$, and $T$ denote the
parallel transport maps respectively in bundles $P$, $\Adp$, and
$\adp$,  over a given path, then we have $S(\phi)=R^{-1}\circ \phi
\circ R$ and $T(x)=dS_e(x)$, the derivative at the identity of $S$
evaluated on the vector $x$.

For $\alpha_i\in\Omega^\ast (M,\adp)$, $1 \leq i \leq n$, define
$\hal(n,i) \in \Omega^\ast (\Delta^n \times LM,ev^\ast \adp)$ by,
\begin{equation*}
\hal(n,i)=T_iev_{n,i}^\ast \alpha_i.
\end{equation*}

Given $\gamma \in LM$ the holonomy along $\gamma$ from $\gamma(0)
$ to $\gamma(1)$ in the principal bundle $P$ gives rise to a
section $hol \in \Gamma (ev^\ast \Adp)$. Note that $\Gamma
(ev^\ast \Adp)$ acts by conjugation on $\Gamma (ev^\ast \Adp)$,
and subsequently on $\Gamma (ev^\ast \adp)$. Now, define
$V^n_{\alpha_1, \cdots, \alpha_n} \in \Omega^\ast  (LM,ev^\ast
\adu)$ as,

\begin{equation}\label{label-wilson-1}
V^0 = 1,
\end{equation}

$$
V^n_{\alpha_1,\cdots,\alpha_n}=\int_{\Delta^n}\hal(n,1)\wedge\cdots
\wedge \hal(n,n), \text{ for } n \geq 1 ,$$ and let,

\begin{equation}\label{formal}
V_{\alpha}=\sum_{n=0}^{\infty} V^n_{\alpha}, ~ \text{where} ~
V^n_{\alpha}=V^n_{\alpha, \cdots, \alpha}.
\end{equation}

It is noteworthy that the above infinite sum is convergent. This
follows from the discussion in Appendix B, based on the simple
fact that the volume of the standard $n$-simplex is $\frac{1}{n!}$
(see Remark \ref{conv}).

Then, consider $W^n_{\alpha_1, \cdots, \alpha_n} \in \Gamma
(ev^\ast \Adp) \times \Omega^\ast  (LM,ev^\ast \adu)$ defined as,
\begin{equation}\label{label-wilson-1}
\begin{split} W^n_{\alpha_1, \cdots, \alpha_n}=(hol, V^n_{\alpha_1, \cdots, \alpha_n}),  \text{ for } n \geq
0,
 \end{split}
\end{equation}
and,
$$ W^n_{\alpha}=(hol, V^n_{\alpha, \cdots, \alpha}) ~ \text{and}
~W_{\alpha}=(hol, V_\alpha).$$

We fix an invariant function $f: G\to \R$ and we shall follow the
notations introduced in (\ref{eq-general-sections-var-func}) and
(\ref{eq-gene-invariant-function}). Note that the maps
$$\hat{f}:\Gamma(\Adp)\otimes \Gamma(\adu)\to
C^{\infty}(M),$$
$$\hat{F}:\Gamma(\Adp)\otimes \Gamma(\adu)\to \Gamma(\adp),$$
naturally induce the maps on the differential forms,
$$\hat{f}:\Gamma(\Adp)\otimes \Omega^\ast (M,\adu)\to \Omega^\ast (M),$$
$$\hat{F}:\Gamma(\Adp)\otimes \Omega^\ast (M,\adu)\to \Omega^\ast (M,\adp),$$
which in turn induce the following maps on the corresponding
pullback bundles,
$$\hat{f}:\Gamma(ev^\ast \Adp)\otimes \Omega^\ast (LM,ev^\ast \adu)\to
\Omega^\ast (LM),$$
$$\hat{F}:\Gamma(ev^\ast \Adp)\otimes
\Omega^\ast (LM,ev^\ast \adu)\to \Omega^\ast (LM, ev^\ast
\adp).$$

\begin{defn}\label{formal2}
For a differential form $\alpha\in \Omega^\ast (M,\adp )$, the
\emph{Wilson loop} $\W_\alpha\in \Omega^\ast (LM,\R)$ is defined
as,
\begin{equation*}\W_\alpha=\hat{f}(W_{\alpha}).
\end{equation*}
\end{defn}

\begin{prop}\label{pro-closed}
If $\alpha \in MC$, then $\W_\alpha\in \Omega^\ast (LM)$ is a
closed form.
\end{prop}
\begin{proof}

Using Stokes theorem  we have,
\begin{equation*}
\begin{split}
d\W_{\alpha}&=d \hat{f}(hol, 1)+ \sum_{n=1}^{\infty} \hat{f}(hol,d
\int_{\Delta^n} \hal(n,1)\wedge\cdots \wedge\hal(n,n))\\&+
\sum_{n=1}^{\infty} \hat{f}( hol, \int_{\partial\Delta^n}
\hal(n,1)\wedge\cdots \wedge\hal(n,n)).
\end{split}
\end{equation*}
Since the connection $\nabla$ is flat,  $\hat f(hol, 1)= f(hol)$
is constant on the connected components of $LM$ therefore $d
\hat{f}(hol)=0$ on $LM$. Actually one can say more, $hol\in
\Gamma(ev^\ast (\Adp))$ is a flat section. Therefore, by Corollary
\ref{app-coro-d-tr}, we have,
\begin{equation*}
\begin{split}
d\W_{\alpha}&=\sum_{n=1}^{\infty} \hat{f}(hol,
\int_{\Delta^n}\dvt(\hal(n,1)\wedge\cdots \wedge\hal(n,n))) \\&+
 \sum_{n=1}^{\infty} \hat{f} ( hol,\int_{\partial\Delta^n} \hal(n,1)\wedge\cdots
\wedge\hal(n,n)).
\end{split}
\end{equation*}
We analyze different terms separately. For the first part we have,
\begin{equation*}
\begin{split}
\dvt(\hal(n,1)\wedge\cdots \wedge\hal(n,n)) &=\sum_{i=1}^n (-1)^i
\hal(n,1)\wedge\cdots \wedge \dvt(\hal(n,i))\wedge \cdots
\wedge\hal(n,n)\\
&=\sum_{i=1}^n (-1)^i \hal(n,1)\wedge\cdots \wedge
(\dt(\alpha))(n,i)\wedge \cdots \wedge\hal(n,n).
\end{split}
\end{equation*}
Hence,
\begin{equation*}
\sum_{n=1}^{\infty}\int_{\Delta^n} \hat{f}(hol,\dvt(
\hal(n,1)\wedge\cdots \wedge\hal(n,n)))=
\sum_{n=1}^{\infty}\sum_{i=1}^{n}\hat{f}( W^n_{\alpha,\cdots,
\dt\alpha, \cdots \alpha}).
\end{equation*}

 As for the second part we have $\partial
\Delta^n=\cup_{i=0}^{n} \Delta^{n-1}_i $, where,
$$\Delta^{n-1}_i= \{(t_0, t_1,\cdots,t_n, t_{n+1}) ~|~ 0=t_0 \leq
t_1\leq \cdots \leq t_i=t_{i+1} \leq\cdots\leq t_n\leq t_{n+1}=1
\},$$ and then for $n\geq1$,
\begin{equation*}
\begin{split}
&
\hat{f}(hol,\int_{\partial\Delta^n}\hal(n,1)\wedge\cdots\wedge\hal(n,n))\\&
=\hat{f}(hol,Ad_{hol}(ev^\ast \alpha)\wedge
\int_{\Delta^{n-1}}\hal(n-1,1)\wedge\cdots \wedge\hal(n-1,n-1))\\&
+\sum_{i=1}^{n-1}(-1)^i \hat{f}(hol, \int_{\Delta^{n-1}}
\hal(n-1,1)\wedge\cdots
\wedge(\alpha\wedge\alpha)(n-1,i)\wedge\cdots \wedge\hal(n-1,n))
\\&+(-1)^n \hat{f}(hol,
\int_{\Delta^{n-1}}\hal(n-1,1)\wedge\cdots \wedge\hal(n-1,n-1)
\wedge ev^\ast \alpha).
\end{split}
\end{equation*}
To calculate the first term in the equality above we have used the
fact that the parallel transport of $\adp$ along a loop $\gamma$
is give by $Ad_{hol}$ where $hol$ is understood as the parallel
transport in the principal bundle $P$. Now, from the invariance of
$f$ and the fact that the components of $\alpha$ have odd degrees,
it follows that the first term and the last term cancel each
other.  Therefore, we have,

\begin{equation*}
\begin{split}
&
\sum_{n=1}^{\infty}\hat{f}(hol,\int_{\partial\Delta^n}\hal(n,1)\wedge\cdots\wedge\hal(n,n))=\\
&\sum_{n=2}^{\infty}\sum_{i=1}^{n-1} (-1)^i\hat{f}(
hol,\int_{\Delta^{n-1}}\hal(n-1,1)\wedge\cdots
\wedge(\alpha\wedge\alpha)(n-1,i)  \wedge \cdots \wedge
\hal(n-1,n-1))\\&=\sum_{n=1}^{\infty}\sum_{i=1}^{n}(-1)^i
\hat{f}(W^n_{\alpha,\cdots, \alpha \wedge \alpha, \cdots
,\alpha}).
\end{split}
\end{equation*}
Then,
\begin{equation*}
\begin{split}
d\W_{\alpha}&= \sum_{n=1}^{\infty}\sum_{i=1}^{n}(-1)^i \hat{f}(
W^n_{\alpha,\cdots, \dt\alpha, \cdots,
\alpha})+\sum_{n=1}^{\infty}\sum_{i=1}^{n}
(-1)^i\hat{f}(W^n_{\alpha,\cdots, \alpha \wedge \alpha, \cdots
,\alpha})\\
&=\sum_{n=1}^{\infty}\sum_{i=1}^{n}(-1)^i\hat{f}(
W^n_{\alpha,\cdots, \dt\alpha+\alpha\wedge\alpha, \cdots
,\alpha}).
\end{split}
\end{equation*}
Note that since $\alpha$ is a sum of forms of odd degrees in
$\Omega^\ast (M,\adu)$, $\alpha\wedge\alpha=1/2[\alpha,\alpha]$.
Therefore,
\begin{equation*}
\begin{split}
d\W_{\alpha}&=\sum_{n=1}^{\infty}\sum_{i=1}^{n}(-1)^i\hat{f}(
W^n_{\alpha,\cdots, \dt\alpha+1/2[\alpha,\alpha], \cdots,
\alpha})=0.
\end{split}
\end{equation*}
\end{proof}

We now explain how $\W_{\alpha}$ represents an equivariant
cohomology class.

\begin{lem} \label{lem-equv-representable}
For  $\alpha\in \MC$, there exists a $\overline{\W}_{\alpha}\in
\cH^\ast (LM)$ such that,
 \begin{equation}\label{eq-W-barW}
 \W_{\alpha}= \mathfrak{e}^\ast  (\w_{\alpha}).
 \end{equation}
 \end{lem}
\begin{proof}
By (\ref {eq-equi-2}) it suffices to show that that
$\W_{\alpha}\in \Omega^\ast (LM) \overset{1\otimes
id}{\hookrightarrow} \Omega^\ast (LM) \otimes
\Omega^\ast(S^\infty)$ is a basic form. Since $\w_{\alpha}$ is a
closed form, it remains to show that $i_v \W_{\alpha}=0$.
 The fact that $ i_v \W_{\alpha}=0$ for
$v=\frac{\partial}{\partial t}$, the fundamental vector field of
the $S^1$-action, is obvious since $ \W_{\alpha}$ is obtained by
pullback of the evaluation maps, $dev_{t_k}(v)=\gamma'(t_k)$ and
$i_{\gamma'(t_k)} i_{\gamma'(t_k)}=0$.  This shows that
$\W_{\alpha}$ determines a cohomology class
$\overline{\W}_{\alpha} \in \cH^\ast (LM)$ which satisfies,
 \begin{equation*}
 \W_{\alpha}= \mathfrak{e}^\ast  (\w_{\alpha}).
 \end{equation*}
\end{proof}
Since the equivariant de Rham cohomology of $LM$ is isomorphic to
the singular one, we can think of $\w_{\alpha}$  as an element in
$H^\ast _{S^1}(LM,\R)$.   For $c\in \H_\ast (LM)$, define
$\Psi_c\in \o(\MC)$ as,
\begin{equation}
\Psi_c(\alpha)=\l c,\w_{\alpha} \r.
\end{equation}
The following proposition implies that $\Psi_c$ descends to
$\o(\Mc)$.
\begin{prop}\label{lem-invariance}
For $c\in H_\ast (LM)$ and $a\in \L_0= \Omega^{2\ast }(M,\adp)$,
\begin{equation*}
L_{\xi_a} \Psi_c=0 \in \o(\MC).
\end{equation*}
\end{prop}
\begin{proof}By Cartan's formula, we have,
\begin{equation*}
\begin{split}
(L_{\xi_a}\Psi_c)(\alpha)&=(i_{\xi_a}d\Psi_c+di_{\xi_a}\Psi_c)(\alpha)=i_{\xi_a}d\Psi_c(\alpha).
\end{split}
\end{equation*}

Let us calculate $d\Psi_c(\alpha)(h)$ for $h\in \L_1=T_\alpha
\L_1$,
\begin{equation*}\begin{split}
d\Psi_c(\alpha)(h)=\frac{d}{dt}\left(\int_c
\w_{\alpha+th}\right)|_{t=0}=\int_c \frac{d}{dt}(
\w_{\alpha+th})|_{t=0}.
\end{split}
\end{equation*}

Note that,

\begin{equation*}
\begin{split}
\frac{d}{dt}\W_{\alpha+th}|_{t=0}=\frac{d}{dt} \sum_{n=1}^{\infty}
\hat{f}(W^n_{\alpha+th,\cdots,\alpha+th})|_{t=0}=
\sum_{n=1}^{\infty}\sum_{i=1}^n \hat{f}(W^n_{\alpha,\cdots,
h,\cdots,\alpha}).
\end{split}
\end{equation*}
So, in order to show that
$L_{\xi_a}\Psi_c(\alpha)=i_{\xi_a}d\Psi_c(\alpha)=0$ we have to
show that $\frac{d}{dt}\W_{\alpha+t h}|_{t=0}$ is an exact form
for $h=\xi_a(\alpha)$,  or in other words,
$$\sum_{n=1}^{\infty}\sum_{i=1}^n
\hat{f}(W^n_{\alpha,\cdots,\xi_a(\alpha),\cdots,\alpha}) ,$$ is exact.
\begin{equation*}
\end{equation*}
Let $Z=\sum_{n=1}^{\infty}\sum_{i=1}^n
(-1)^i\hat{f}(W^n_{\alpha,\cdots,a,\cdots,\alpha}) \in \Omega^\ast
(LM,\R)$. Then, by a calculation similar to that presentation in
the proof of Proposition \ref{pro-closed} we show that for $a\in
\L_+$, \begin{equation*}
\begin{split}
dZ=&\sum_{n=1}^{\infty}\sum_{i=1}^n
(-1)^{2i}(\hat{f}(W^n_{\alpha,\cdots, \tiny\underbrace{\dt a}_{i}
,\cdots,\alpha})+\hat{f}(W^n_{\alpha,\cdots,\tiny\underbrace{
\alpha \wedge a}_{i},
 \cdots,\alpha})-\hat{f}(W^n_{\alpha,\cdots, \tiny\underbrace{a \wedge \alpha}_{i},
 \cdots,\alpha})+\\&\sum_{n=1}^{\infty}\sum_{i=1}^n (-1)^i[\sum_{j<i}
 (-1)^j(\hat{f}(W^n_{\alpha,\cdots,\tiny\underbrace{ \dt \alpha}_{j}, \cdots\tiny\underbrace{,a,}_{i}
 \cdots,\alpha})+\hat{f}(W^n_{\alpha,\cdots, \tiny\underbrace{\alpha\wedge
 \alpha}_{j},
 \cdots\tiny\underbrace{,a,}_{i}
 \cdots,\alpha}))+\\&\sum_{i<j\leq
 n}(-1)^j(\hat{f}(W^n_{\alpha,\cdots
\tiny\underbrace{,a,}_{i}\cdots ,\tiny\underbrace{\dt \alpha}_{j},
 \cdots,\alpha}))]+\hat{f}(W^n_{\alpha,\cdots\tiny\underbrace{,a,}_{i}\cdots,\tiny\underbrace{\alpha \wedge
 \alpha}_{j},
 \cdots,\alpha}))
 \\&=\sum_{n=1}^{\infty}\sum_{i=1}^n (-1)^{2i} (\hat{f}(W^n_{\alpha,\cdots,
\tiny\underbrace{\dt
a}_{i},\cdots,\alpha})+\hat{f}(W^n_{\alpha,\cdots,
\tiny\underbrace{\alpha \wedge a}_{i},
 \cdots,\alpha})-\hat{f}(W^n_{\alpha,\cdots, \tiny\underbrace{a \wedge
 \alpha}_{i},
 \cdots,\alpha}))+\\&\sum_{n=1}^{\infty}\sum_{i=1}^n(-1)^i (\sum_{j<i}
  (-1)^j\hat{f}(W^n_{\alpha,\cdots, \tiny\underbrace{\dt \alpha+1/2[\alpha,\alpha]}_{j}, \cdots\tiny\underbrace{,a,}_{i}
 \cdots,\alpha})\\&+\sum_{i<j\leq n}(-1)^j\hat{f}(W^n_{\alpha,\cdots
 \tiny\underbrace{,a,}_{i}\cdots, {\tiny\underbrace{\dt \alpha+1/2[\alpha,\alpha]}_{j}},
 \cdots,\alpha})).
\end{split}
\end{equation*}

As $\dt \alpha+\frac{1}{2}[\alpha,\alpha]=0$ for $\alpha\in \MC$,
therefore,
\begin{equation*}
\begin{split}
dZ=&\sum_{n=1}^{\infty}\sum_{i=1}^n \hat{f}(W^n_{\alpha,\cdots,
\dt a +[\alpha,a],\cdots,\alpha})\\&
=\sum_{n=1}^{\infty}\sum_{i=1}^n \hat{f}(W^n_{\alpha,\cdots,
\xi_a(\alpha),\cdots,\alpha})\\&
=\frac{d}{dt}\W_{\alpha+t\xi_a(\alpha)}|_{t=0}.
\end{split}
\end{equation*}
This proves the claim.
\end{proof}
By Proposition \ref{lem-invariance}, $\Psi_a:MC \to \mathbb R$ is
constant along the orbits of $\G$.
\begin{defn} \label{Psi}
The generalized holonomy map  $\Psi: \H_{2\ast }(LM) \to \o (\Mc)$ is the map $ a \mapsto
\Psi_a$, where,
$$ \Psi_a(\alpha)=\langle a, \w_{\alpha}\rangle.$$
\end{defn}

\section{Hamiltonian vectors field and Poisson bracket}\label{Hamilt}

Let $\Psi_a$, and $\Psi_b$ denote the holonomy functions
associated to the equivariant homology classes $a, b\in \H_{2\ast
}(LM)$, and let $X_a$ and $X_b$ denote their corresponding
Hamiltonian vector fields. We calculate the Poisson bracket of
holonomy functions by evaluating the symplectic form on their
corresponding Hamiltonian vector fields.  Note that the Zarisky
tangent space of $\MC/\G$ at a class represented by $\alpha$ can
be identified with the cohomology whose differential is $\nabla+
\alpha\wedge\cdot$ or in other words $\adp$ is equipped with the
flat connection $\nabla+ \alpha \wedge\cdot$. The later is a
differential as $\alpha$ satisfies the Maurer-Cartan equation. We
denote the cohomology of this differential by $H^\ast
(M,\adp_{\alpha})$, and similarly the corresponding vector bundle
homology by $H_\ast (M,\adp_{\alpha})$.

We borrow the following lemma from \cite{G2} to obtain a formula for the Hamiltonian
vector fields.

\begin{lem}\label{lem-PD-iso} The following diagram is
commutative,
\begin{equation}
\xymatrix{H^k(M,\adp_{\alpha}) \ar^-{\tilde{\omega}} [d]\ar^-{\tilde{\theta}} [rd]\ar^-{\cap [M]}[r] & H_{m-k}(M,\adp_{\alpha})\ar^-{\tilde\eta} [d]\\
H^{m-k}(M,\adp_{\alpha})^\ast  & H^{m-k}(M,\coad _{\alpha})^\ast
\ar^-{\tilde{b}^t }[l] }
\end{equation}
\end{lem}
where,
\begin{romlist}
\item $ [M]\in H_d(M)$ is the fundamental class of $M$, \item
$\tilde{\omega}$ is induced by $\omega:
H^k(M,\adp_{\alpha})\otimes
H^{m-k}(M,\adp_{\alpha})\overset{\int_M  \l\cdot,\cdot\r} {
\longrightarrow} \R$, \item $\tilde{\theta}$ is induced by
$\theta:H^k(M,\adp_{\alpha})\otimes
H^{m-k}(M,\coad_{\alpha})\overset{\small{(\int_M \cdot \wedge
\cdot) \otimes \text{eval}}}{\longrightarrow} \R$, \item
$\tilde{\eta}$ is induced by $\eta:H_{m-k}(M,\adp_{\alpha})\otimes
H^{m-k}(M,\coad_{\alpha})\overset{\text{eval}}{\longrightarrow}
\R$, \item $\tilde{b}^t$ is the transpose of the map
$\tilde{b}:H^{m-k}(M,\adp_{\alpha})\to H^{m-k}(M,\coad_{\alpha})$
induced by $\l\cdot, \cdot\r:\g \times \g \to \R$.
\end{romlist}
Here, $\coad_\alpha $ is the bundle associated to the
representation $Ad^\ast : G \to Aut(\mathfrak g^\ast )$ induced by
the conjugation action, with the differential induced by $\nabla +
\alpha \wedge \cdot$.

\begin{lem}\label{lem-vam-vector}The Hamiltonian vector field $X_a(\alpha) \in H^\ast (M, \adp_{\alpha})$, for $a \in H_\ast (LM)$, satisfies the following equation
\begin{equation*}
\begin{split}
PD(X_a(\alpha))=ev_\ast  (\mathfrak{m}_\ast a\cap
\hat{F}(W_{\alpha}))=\sum_{k=0}^{\infty} ev_\ast
\left(\mathfrak{m}_\ast a\cap \hat{F}(W^k_{\alpha})\right),
\end{split}
\end{equation*}
where $PD$ denotes the Poincar\'e duality map $H^\ast
(M,\adp_{\alpha}) \overset{\tiny \cdot \cap  [M]}{\longrightarrow}
H_\ast (M,\adp_{\alpha})$.
\end{lem}
\begin{proof} For the Hamiltonian vector field $X_a(\alpha)$ and $h \in T_{\alpha}\Mc$,
\begin{equation*}
\begin{split}
\omega(X_a(\alpha),h)=&\frac{d}{dt}\Psi_a(\alpha+th)|_{t=0}=\frac{d}{dt}
\l a,\overline{\W}_{\alpha+th}\r|_{t=0}= \l a, \frac{d}{dt}
,\overline{\W}_{\alpha+th}\r |_{t=0}.
\end{split}
\end{equation*}
Note that $ \w_{\alpha+th}$ is represented by the differential
form,
$$\W_{\alpha+th}=\sum_{k=0}^{\infty}\hat{f}(W^k_{\alpha+th,\cdots
,\alpha+th}).$$ Because of the flatness of $\nabla$,
$d\hat{f}(hol)=0$, therefore,
 \begin{equation*}
\begin{split}
\mathfrak{e}^\ast (\frac{d}{dt}
\overline{\W}_{\alpha+th}|_{t=0})=\sum_{k=1}^{\infty}\sum_{i=1}^{n}\hat{f}(W^k_{\alpha\cdots\tiny\underbrace{,h,}_{i}\cdots,\alpha}).
\end{split}
\end{equation*}
which by Theorem 2.1 \cite{GJP} is,
\begin{equation*}
\begin{split}
\sum_{k=1}^{\infty}\Delta^\ast  \hat{f}(hol, V^k_{\alpha}\wedge
ev^\ast h)=\sum_{k=0}^{\infty}\mathfrak{e}^\ast \mathfrak{m}^\ast
\hat{f}(hol, V^k_{\alpha}\wedge ev^\ast h).
\end{split}
\end{equation*}

Hence by (\ref{eq-general-variation}) and (\ref{eq-gene-invariant-function}),
\begin{equation*}
\begin{split}
\omega(X_a(\alpha),h)&=\l a, \mathfrak{m}^\ast \sum_{k=0}^{\infty}
\hat{f}(hol, V^k_{\alpha}\wedge ev^\ast h) \r \\ &=\l
\mathfrak{m}_\ast a, \sum_{k=0}^{\infty} \hat{f}(hol,
V^k_{\alpha}\wedge ev^\ast h) \r\\ &= \l \mathfrak{m}_\ast a\cap
\sum_{k=0}^{\infty} \hat{F}(W^k_{\alpha}) , ev^\ast  h\r
\\&=\l ev_\ast (\mathfrak{m}_\ast a \cap \sum_{k=0}^{\infty} \hat{F}(W^k_{\alpha})), h
\r.
\end{split}
\end{equation*}
Therefore, by Lemma \ref{lem-PD-iso},
\begin{equation*}\begin{split}
ev_\ast (\mathfrak{m}_\ast a\cap \sum_{k=0}^{\infty}
\hat{F}(W^k_{\alpha}))
&=\eta^{-1}\circ (\tilde{b}^t)^{-1}(d_{\alpha}\Psi_a)\\
&=\tilde{\omega}^{-1}(d_{\alpha}\Psi_a)\cap [M]\\ & \\&=
X_a(\alpha)\cap [M].
\end{split}
\end{equation*}

\end{proof}

The following lemma is a version of the multiplicative property of
the formal power series parallel transport (see \cite{Ch1, Ch,
M}). See (\ref{star-operator}) to recall the definition of the
$\ast$ product.
\begin{lem}\label{lem-pull-back-holonomy} Let  $\varepsilon: LM\times_M LM\to LM$ be the concatenation map,
then, $$\varepsilon^\ast (W_{\alpha})= pr_1^\ast W_{\alpha}\ast
pr_2^\ast W_{\alpha},$$ where $pr_{1}$ and $pr_{2}$ are the
projections on the first and the second factors,
\begin{equation*} \xymatrix{& LM\times_M LM\ar^-{pr_1}[dl]\ar_-{pr_2}[dr]\\LM&&LM
}
\end{equation*}
\end{lem}
\begin{proof}
The proof of Proposition 1.5.1 in \cite{Ch1} can be adopted to
prove that that for two loops $\gamma_1$ and $\gamma_2$ with
identical marked point,

\begin{equation} \label{first}
V^n_{\alpha}(\gamma_1\circ\gamma_2)=\sum_{i+j=n}
Ad_{hol(\gamma_2)}V^i_{\alpha}(\gamma_1)\wedge
V^j_{\alpha}(\gamma_2).
\end{equation}
Also $hol(\gamma_1\circ \gamma_2)=hol(\gamma_2)\circ
hol(\gamma_1)$. To see the identity (\ref{first}) more clearly,
recall that the parallel transport in $\adp$ along $\gamma_2$ is
given by $Ad_{hol(\gamma_2)}$, and that  in computing the
differential form $V^n_{\alpha}(\gamma_1\circ\gamma_2)$ (see
section \ref{GenHol}) the parallel transport $T_i$  along
$\gamma_1\circ \gamma_2$ acts on all the pullbacks
$ev_{t_i}^*\alpha$. Split the terms $\hal(n,i)$ into two
collections according to whether they are obtained by the pullback
with respect to $ev_t$ for $0\leq t \leq 1/2$ or $1/2\leq t \leq
1$. Then, note that one needs to apply an extra action of
$Ad_{hol(\gamma_2)}$ on the first collection to account for the
parallel transport along the second half of $\gamma_1\circ
\gamma_2$, that is $\gamma_2$.
\end{proof}

Recall Definition \ref{Psi} wherein the generalized holonomy map $\Psi: \H_{2\ast }(LM) \to \o (\Mc)$ is described as the map $ a \mapsto
\Psi_a$, where, $ \Psi_a(\alpha)=\langle a, \w_{\alpha}\rangle$. For a version of the following theorem in the BRST setup see \cite{CFP}.

\begin{thm}\label{main-thm} For $G=\Gl (n,\C)$ or $\Gl (n,\R)$, the generalized holonomy
map,
$$\Psi:(\H_{2\ast }(LM),[\cdot,\cdot]) \to
(\o(\Mc),\{\cdot,\cdot\}),$$ is a map of Lie algebras.
\end{thm}

\begin{proof}Suppose that  $G=\Gl (n,\C)$, it is exactly the same argument for $G=\Gl (n,\R)$. We prove that for $a,b\in \H_{2\ast }(LM)$ and $\alpha \in \Mc$,
$$
\Psi_{[a,b]}(\alpha)=\{\Psi_a,
\Psi_b\}(\alpha)=\omega(X_a(\alpha),X_b(\alpha)).
$$
The left hand side is,
\begin{equation*}
\begin{split}
 \Psi_{[a,b]}(\alpha)&=\Psi_{(\mathfrak{e}_\ast (\mathfrak{m}_\ast a\bullet \mathfrak{m}_\ast b)}(\alpha)\\ &=\l \mathfrak{e}_\ast (\mathfrak{m}_\ast a\bullet \mathfrak{m}_\ast b),\w_{\alpha} \r
 \\&=\l \mathfrak{m}_\ast a\bullet \mathfrak{m}_\ast b,\mathfrak{e}^\ast (\w_{\alpha})\r\\  &=\l \mathfrak{m}_\ast a\bullet \mathfrak{m}_\ast b,\W_{\alpha}\r \\
 &=\l \varepsilon_\ast \circ \tau_\ast (\mathfrak{m}_\ast a\times \mathfrak{m}_\ast b), \hat \tr(W_{\alpha})\r\\
 & =\l \mathfrak{m}_\ast a\times \mathfrak{m}_\ast b, \hat \tr( \tau^\ast  \circ \varepsilon^\ast
 (W_{\alpha}))\r,
\end{split}
\end{equation*}
which using Lemma \ref{lem-pull-back-holonomy} equals,
\begin{equation*}
\begin{split}
\int_{ \mathfrak{m}_\ast a\times \mathfrak{m}_\ast b}( ev\times
ev)^\ast (U) \wedge \hat{\tr}( pr_1^*W_{\alpha}\ast
pr_2^*W_{\alpha}),
\end{split}
\end{equation*}
where $U\in \Omega^\ast (M\times M)$ is a differential form
supported in a neighborhood of the diagonal realizing the Thom
isomorphism. Now, by Lemma \ref{lem-vam-vector},
\begin{equation*}
\begin{split}
\{\Psi_a,\Psi_b\}(\alpha)&=\omega(X_a(\alpha),X_b(\alpha))
\\ &=\int_M
\tr(X_a(\alpha)\wedge X_b(\alpha))\\
&=\int_M
 \tr( PD^{-1} (ev_\ast (\mathfrak{m}_\ast a\cap \hat{F}(W_{\alpha}))\wedge  PD^{-1}(ev_\ast (\mathfrak{m}_\ast a\cap
 \hat{F}(W_{\alpha}))).
 \end{split}
\end{equation*}
Using Remark \ref{remark-important} and the inclusion $\Omega^\ast
(LM, ev^\ast \Adp)\hookrightarrow \Omega^\ast (LM,ev^\ast \adp)$
induced by $\Gl(n,\C)\hookrightarrow \g l(n,\C)$, and by Corollary
\ref{cor-variantion-GL},
 $$
 \hat{F}(W_{\alpha})=hol \cdot \Pi V_{\alpha}\in \Omega^\ast
 (LM,ev^\ast \adp).
 $$
Thus, the above Poisson bracket equals,
 \begin{equation*}
\begin{split}
 &\int_{\tr(ev_\ast (\mathfrak{m}_\ast a\cap \hat{F}(W_{\alpha})) \times ev_\ast (\mathfrak{m}_\ast a\cap \hat{F}(W_{\alpha})))} U=\int_{(\mathfrak{m}_\ast a\times \mathfrak{m}_\ast b)\cap \tr (\hat{F}(W_{\alpha})
 \times \hat{F}(W_{\alpha}))} (ev\times ev)^\ast U\\&=\int_{\mathfrak{m}_\ast a\times \mathfrak{m}_\ast b} (ev \times ev)^\ast U \wedge \tr (pr_1^\ast \hat{F}(W_{\alpha})\wedge pr_2^\ast \hat{F}(W_{\alpha}))\\&=
 \int_{\mathfrak{m}_\ast a\times \mathfrak{m}_\ast b} (ev \times ev)^\ast U \wedge \tr (pr_1^\ast (hol \cdot \Pi V_{\alpha})\wedge pr_2^\ast (hol \cdot \Pi V_{\alpha}))\\&=
 \int_{\mathfrak{m}_\ast a\times \mathfrak{m}_\ast b} (ev \times ev)^\ast U \wedge \tr( (pr_1^\ast hol\cdot pr_1^\ast \Pi V_{\alpha})\wedge
 (pr_2^\ast hol\cdot pr_2^\ast \Pi V_{\alpha})).
 \end{split}
\end{equation*}

But, it follows immediately from Remark \ref{remark-important} and
line (\ref{calcul-trace}) that,
$$\tr ((pr_1^\ast hol\cdot pr_1^\ast \Pi V_{\alpha})\wedge (pr_2^\ast hol\cdot pr_2^\ast \Pi
V_{\alpha}))= \hat{\tr} (pr_1^\ast W_{\alpha}\ast pr_2^\ast
W_{\alpha}).$$ This proves the theorem.
\end{proof}

\section{String bracket for unoriented strings}\label{unor}
In this section we calculate the Poisson bracket for the
generalized holonomy associated to  a principal $G$-bundle
equipped with a flat connection, where $G$ is one of the proper
reductive subgroups of  $\Gl(n,\C)$ such as $\O(n,\R)$,
$\O(n,\C)$, $\mathrm{U}(n)$, $\Sl(n,\R)$, or $\Sl(n,\C)$. Let
$i:LM \to LM $ be the involution which reverses the orientation of
a loop. We continue to denote the induced map on the homology
$H_\ast (LM)$  and the equivariant homology $\H_\ast (LM)$ by
$a\mapsto i(a)$. We use $i^\ast $ for the induced map on the
cohomologies. It is a direct check that the latter maps commute
with the maps $\mathfrak{e}_\ast $ and $\mathfrak{m}_\ast $
introduced in section \ref{StringBracket}.

\begin{thm}\label{thm-main-UN}
For $G= \O(p, q)$, $\O(n, \C)$,$\U(p,q)$, $\Sp(n,\R)$ and
$\Sp(p,q)$,
\begin{equation}
\{\Psi_a,\Psi_b\}=\frac{1}{2}(\Psi_{ [a,b]}-\Psi_{ [a,i(b)]}),
\end{equation}
where $a,b\in \H_{2\ast }(LM)$.
\end{thm}
\begin{proof}
Consider
$$
 V_{\alpha}^{t}=\sum_{n=0}^{\infty}(-1)^n \int_{\Delta^n}\hal(n,n)\wedge\cdots \wedge
 \hal(n,1).
$$
Notice that first term of the above summation is $1$ (compare (\ref{label-wilson-1})). In a manner similar to the proof
of Theorem \ref{main-thm},
\begin{equation*}
\begin{split}
\{\Psi_a, \Psi_b\}(\alpha)&=\int_{\mathfrak{m}_\ast a\times
\mathfrak{m}_\ast b} (ev\times ev)^\ast U \wedge \tr (pr_1^\ast
\hat{F}(W_{\alpha})\wedge pr_2^\ast \hat{F}(W_{\alpha}))),
\end{split}
\end{equation*}

By Corollary \ref{variation-UN-envelop},
$$
\hat{F}(W_{\alpha})=\frac{1}{2}(hol\cdot \Pi V_{\alpha}-\Pi
V_{\alpha}^{t} \cdot hol^{-1}) \in \Omega^\ast (LM,ev^\ast \adp)
,$$ and all the terms above make sense using the  fact  $G$ is one
of the subgroups listed in the theorem.  Since for every $g\in G$,
$\Re\tr (A)=\Re \tr (A^{-1})$ it follows,

\begin{equation*}
\begin{split}
\{\Psi_a, \Psi_b\}(\alpha)&=\frac {1}{2}\int_{\mathfrak{m}_\ast
a\times \mathfrak{m}_\ast b} (ev\times ev)^\ast U \wedge \tr
(pr_1^\ast (hol) pr_1^\ast (\Pi V_{\alpha})\wedge pr_2^\ast
(hol)pr_2^\ast (\Pi V_{\alpha}))\\&-\frac
{1}{2}\int_{\mathfrak{m}_\ast a\times \mathfrak{m}_\ast b}
(ev\times ev)^\ast U \wedge \tr (pr_1^\ast (hol) pr_1^\ast (\Pi
V_{\alpha})\wedge pr_2^\ast (\Pi V_{\alpha}^{t})pr_2^\ast
(hol^{-1})).
\end{split}
\end{equation*}
Also, as in the proof of Theorem \ref{main-thm},

\begin{equation*}
\begin{split}
 \frac{1}{2}(\Psi_{[a,b ]}(\alpha)-\Psi_{[a,i(b) ]}(\alpha)) &=\frac{1}{2}\int_{\mathfrak{m}_\ast a\times \mathfrak{m}_\ast b} (ev\times ev)^\ast (U) \wedge \tr (pr_1^\ast W_{\alpha}\ast pr_2^\ast W_{\alpha})\\
 & -\frac{1}{2}\int_{\mathfrak{m}_\ast a\times \mathfrak{m}_\ast i(b) } (ev\times ev)^\ast (U) \wedge \tr (pr_1^\ast W_{\alpha}\ast pr_2^\ast
 W_{\alpha}).
\end{split}
\end{equation*}

Similarly to the proof of Theorem \ref{main-thm},
\begin{equation*}\begin{split}
&\int_{\mathfrak{m}_\ast a\times \mathfrak{m}_\ast b}\tr
(pr_1^\ast (hol) pr_1^\ast \Pi V_{\alpha}\wedge pr_2^\ast
(hol)pr_2^\ast \Pi V_{\alpha})\\&=\int_{\mathfrak{m}_\ast a\times
\mathfrak{m}_\ast b} (ev\times ev)^\ast (U) \wedge \hat{\tr}
(pr_1^\ast W_{\alpha}\ast pr_2^\ast W_{\alpha}).
\end{split}
\end{equation*}
So, we only have to show that,
\begin{equation}\label{petit-claim-unitary}
\begin{split}
&\int_{\mathfrak{m}_\ast a\times \mathfrak{m}_\ast  b}\tr
(pr_1^\ast (hol) pr_1^\ast \Pi V_{\alpha}\wedge
pr_2^\ast \Pi V^t_{\alpha}pr_2^\ast (hol^{-1}))\\
&=\int_{\mathfrak{m}_\ast a\times \mathfrak{m}_\ast i(b)}(ev\times
ev)^\ast (U) \wedge \hat{\tr}(pr_1^\ast W_{\alpha}\ast pr_2^\ast
W_{\alpha}).
\end{split}
\end{equation}
Note that $ev \circ i=ev$, and therefore $(ev\times ev)^\ast U
=(id \times i)^\ast (ev\times ev)^\ast U$. Also,
\begin{equation*}
\begin{split}
&{\int_{\mathfrak{m}_\ast a\times \mathfrak{m}_\ast i(b)}} (ev\times ev)^\ast U \wedge \hat{\tr} (pr_1^\ast W_{\alpha}\ast pr_2^\ast W_{\alpha})\\
&={\int_{\mathfrak{m}_\ast a\times \mathfrak{m}_\ast b}} (id\times i)^\ast ((ev\times ev)^\ast (U) \wedge \hat{\tr} (pr_1^\ast W_{\alpha}\ast pr_2^\ast W_{\alpha}))\\
&={\int_{\mathfrak{m}_\ast a\times \mathfrak{m}_\ast b}} (ev\times
ev)^\ast U \wedge \hat{\tr} (pr_1^\ast W_{\alpha}\ast pr_2^\ast
(i^\ast (W_{\alpha})).
\end{split}
\end{equation*}
Note\footnote{We remind the reader that in our convention
$Ad_g(h)=g^{-1}hg$.} that the pullback $i^\ast
(W_{\alpha})=(hol^{-1},Ad_{hol^{-1}}V^t_\alpha)$, and therefore by
Remark \ref{remark-important} and equation (\ref{calcul-trace}),
\begin{equation*}
\begin{split}
 \hat{\tr} (pr_1^\ast W_{\alpha}\ast pr_2^\ast (i^\ast (W_{\alpha})))&=\hat{\tr}(pr_1^\ast (hol,V_{\alpha})\ast pr_2^\ast (hol^{-1}, Ad_{hol^{-1}}V^t))\\
 &=\tr(pr_1^\ast (hol) \cdot pr_1^\ast \Pi V_{\alpha}\wedge pr_2^\ast \Pi V^t_{\alpha}\cdot pr_2^\ast
 (hol^{-1})).
\end{split}
\end{equation*}
This verifies equation (\ref{petit-claim-unitary}).
\end{proof}

For a homology class $a\in \H_\ast (LM)$, let $\bar{a}=a+i(a)$
which could be thought of as a homology class of unoriented loops.

\begin{lem}\label{lem-unoriented-map} For $a\in \H_\ast (LM)$,  $a\mapsto i(a)$ is a Lie algebra map.
\end{lem}
\begin{proof}
For $a,b\in H_\ast (LM)$, $i(a)\bullet  i(b)= \varepsilon\circ
\tau( i(a)\times i(b))=\varepsilon\circ (i,i)\circ \tau(a\times
b)$ as the Thom collapsing map commutes with the orientation
reversing map. Also,
\begin{equation}
\varepsilon\circ (i,i)\circ \tau(a\times b)=i \circ \varepsilon
\circ \tau(b\times a)=(-1)^{|a||b|} b\bullet a.
\end{equation}
As for the Lie bracket,

\begin{equation*}
\begin{split}
[i(a),i(b)]&=(-1)^{|a|} \mathfrak{e}_\ast (\mathfrak{m}_\ast
(i(a))\bullet \mathfrak{m}_\ast (i(b)))\\ & =(-1)^{|a|}
(-1)^{|a||b|}\mathfrak{e}_\ast (i(\mathfrak{m}_\ast b\bullet
\mathfrak{m}_\ast a))
\\ &=(-1)^{|a|} (-1)^{|a||b|} i(\mathfrak{e}_\ast (\mathfrak{m}_\ast b\bullet \mathfrak{m}_\ast a))\\ &=(-1)^{|a|+|a||b|+|b|} i([b,a])\\
&=(-1)^{|a|+|a||b|+|b|} (-1)^{|b|+|a|+|a||b|} i([a,b ])\\
&=i([a,b]).
\end{split}
\end{equation*}

\end{proof}

Then, it follows from  the lemma above that,

\begin{equation*}
\begin{split}
 [\bar{a},\bar{b}]&= [a,b]+[i(a),i(b)]+  [i(a),b]+[a,i(b)]\\&=[a,b]+i([a,b])+  i([a,i(b)])+[a,i(b)]\\
 &=\overline{[a,b]}+\overline{[a,i(b)]}\\&=\overline{[a,b]}+\overline{[a,i(b)]}.
\end{split}
\end{equation*}
This  motivates a new Lie bracket on $\H_\ast (LM)$ defined by,
 \begin{equation}
[a,b]_1=\frac{1}{2}\overline{[\bar{a},\bar{b}]}=\frac{1}{2}([a,b]+[a,i(b)]).
\end{equation}
Now, using the above observation, Theorem \ref{thm-main-UN} can be
reformulated as follows:

\begin{thm}\label{thm1-main-UN}
For $G= \O(p, q)$, $\O(n, \C)$,$\U(p,q)$, $\Sp(n,\R)$ and
$\Sp(p,q)$,
\begin{equation}
\{\Psi_{a},\Psi_{b}\}=\Psi_{ [a,b]_1},
\end{equation}
In other words, the map $\H_{2\ast }(LM) \to \o(\Mc)$, defined by
$a\mapsto \Psi_{a}$, is a map of Lie algebras.
\end{thm}

\section{The $2$-dimensional case}\label{2dim}

In this section we will show how our construction relates to the
original construction of Goldman's in \cite{G1, G2} where he
studied the space of flat connection from a representation theory
viewpoint.

Let $M$ be a surface with fundamental group $\pi=\pi_1(M)$. Then
the representation variety $Rep(\pi,G)$, as a set, is in a
one-to-one correspondence with the isomorphism classes of flat
$G$-bundles over $M$ (see A.1).  It is  known that $Rep(M,\pi)$ is
a symplectic space \cite{G1,AB}. A connected component of
$Rep(\pi,G)$ determines an isomorphism class of principal bundles
$G \to P\to M$, and the elements of the of each connected
component are in a one-to-one correspondence with the flat
connections on the corresponding bundle modulo the gauge group.

 Let us fix a $G$-bundle $P$, or equivalently, a connected component of $Rep(\pi,G)$}. Then the space of all connections on $G \to P \to M$ is an affine space modeled over the vector space $\Omega^1(M, \adp)$. In particular, for a fixed
connection $\nabla_0$, the map $\alpha \mapsto \nabla_0 + \alpha$
form $\Omega^1$ to the space of all connections is a bijection.
Assuming that $\nabla_0$ is a flat connection, this map gives a
one-to-one correspondence between the set of those $1$-forms which
satisfy the Maurer-Cartan equation and that of all flat
connection. In fact, the above map induces a bijection between the
equivalence classes of the solution of the Maurer-Cartan equation
in $\Omega^1(M, \adp)$, denoted by $\Mc_1$, and the space of flat
connections modulo gauge equivalence. In fact, if $(M, \omega)$ is
a symplectic manifold of dimension $m=2d$ for which Hard Lefschetz
theorem holds, then the moduli space of flat connections on $P$
may be viewed as a symplectic substack of $\Mc$ (see section 1.7
of \cite{GG}). More specifically, the map $\nabla_0+\alpha \mapsto
1/2(\alpha + \alpha\wedge\omega^{d-1})$ identifies the moduli
space of flat connections with a symplectic sub-stack of $\Mc$
(compare \cite{Kar}).

 Now, assume $M$ be an
orientable manifold of dimension $2$. Therefore, for $i \geq 3$,
we have $\Omega^i(M, \adp)=0$. We therefore have $\Mc_1=\Mc$.  Moreover  the symplectic form considered by Goldman for the representation variety corresponds  to the one consider in section 4 (see \cite{G1,AB}).

So to make our  point it only remains to show how the restriction
of the map $\Psi: \H_{2\ast } (L M)\to \mathcal{O} (\Mc)$ to the
equivariant $0^{th}$-homology $\H_0(LM)$ is precisely the map
discovered by Goldman in \cite{G2}. First thing to recall is that
$\H_0(LM)$ equals $\mathbb R \hat \pi$, the vector space generated
by the set of all conjugacy classes of the fundamental group
$\pi_1(M)$. Now, for an element $a \in \H_0(LM)$ look at $\Psi_a
\in \mathcal{O} (\Mc)$, where $\Psi_a(\alpha)=\langle a,
\w_{\alpha}\rangle$. Recall that $\W_\alpha=\hat{f}(W_{\alpha})$,
$W_{\alpha}=(hol, V_\alpha)$, and $V_{\alpha}=\sum_{n=0}^{\infty}
V^n_{\alpha}, ~ \text{where} ~ V^n_{\alpha}=V^n_{\alpha, \cdots,
\alpha}$. Since $\alpha$ has to be a $1$-form, the discussion of
Appendix \ref{Appen-Hol} on iterated integrals says that for a
loop in $\gamma \in LM$, $ \Psi_a(\gamma)=\langle \gamma,\w_\alpha
\rangle$ is nothing but the invariant function $f$ applied to the
holonomy of the connection $\nabla_0 +\alpha$ along the loop
$\gamma$. Notice that, via the one-to-one correspondence of
Proposition \ref{flatconn-rep}, this is exactly the function
$f_{[\gamma]}$ on $Hom(\pi , G)/G$ where $[\gamma]\in \hat{\pi}$
is the free homotopy class determined by the curve $\gamma$. We
therefore have,

\begin{prop} For a closed orientable surface $M$, the restriction
of the Lie algebra map $\Psi: H_{2\ast}^{S^1}(LM) \to \mathcal
O(\Mc)$, to the $0^{th}$-equivariant homology $H_0^{S^1}(LM)$ is
the same is the map $\gamma \mapsto f_{[\gamma]}$, from $\R \hat
\pi$ to the Poisson algebra of functions on the symplectic space
$Hom(\pi, G)/G$, discovered by Goldman in \cite{G2}.
\end{prop}

\appendix

\section{Vector bundles and flat connections} \label{Vect}
In this appendix first we review some basic facts relating  the
space of representations of $\pi_1(M)$ the fundamental group of a
manifold $M$ in a Lie group $G$ and the space of the isomorphism
classes of principal bundles on $M$ equipped with a flat
connection. Second we prove some basic facts which at end enables
us to show that $\Omega^\ast (M,\adp)$ is an example of DGLA as it
is defined in Section \ref{section-MC}. In the end we recall some
homology and cohomology with coefficient in a flat bundle. For
more details, we refer the reader to \cite{Gr}. Notice that a
connection on a principal bundle $G\to P\to M$ gives rise to a
connection on the associated vector bundles. Therefore, we denote
the connections on its associated vector bundles with the same
notation $\nabla$. The covariant derivative $\nabla: \Omega^0(M,
E)\to \Omega^1(M,E)$ can be extended by Leibnitz rule to a
differential $d_{\nabla}:\Omega^{*}(M,E)\to \Omega^{*+1}(M,E)$.
The flatness of the connection translate into the equation
$d_{\nabla}^2=0$.

Let $\pi_1(M)$ denote the fundamental group of a manifold $M$ with
a based point $b$. Let $Hom(\pi_1(M), G)$ denote the set of all
group homomorphisms from $\pi_1(M)$ to a Lie group $G$. Let
$\mathcal F$ denote the set of all pairs $(P, \nabla)$ of
principal $G$-bundles $P$ equipped with a flat connection $\nabla$
up to bundle isomorphisms. The universal covering $\pi_1(M) \to
\tilde{M} \to M$ may be viewed as a principal $\pi_1(M)$-bundle
over $M$ and therefore given a homomorphism $\rho: \pi_1(M) \to G$
one can construct the associated bundle $\tilde{M}\times_\rho G$.
The canonical connection given by unique lifting property will
induce a connection on $\tilde{M}\times_\rho G$. This establishes
a map $\Phi: Hom(\pi_1(M), G) \to \mathcal F$ by $\Phi(\rho)=
\tilde{M}\times_\rho G$. The following statement is the Theorem
6.60 of \cite{Mo}.
\begin{prop} \label{flatconn-rep} The map that associates to any flat $G$-bundle over $M$, its holonomy is a one-to-one correspondence between
the set of conjugacy classes of flat $G$-bundles and the set of
conjugacy classes of homomorphisms $\rho : \pi_1 \to G$.
\end{prop}

In this paper, we are mainly interested in the adjoint bundle $\g
\to \adp\to M$, and the universal enveloping bundle $\u \to
\adu\to M$, associated with the adjoint representation $Ad:G\to
Aut(\g)$ and $Ad_u:G\to Aut(\u)$. Below are a few applications of
the previous results.

\begin{cor}\label{App-compatible-bracket}
For $x,y\in  \Omega^\ast (M,\adp)$ we have,
$$
\dt [x,y]= [\dt x,y ]+ (-1)^{|x|} [x,\dt y].
$$
\end{cor}
\begin{proof}
The claim may be verified locally.  In a local trivialization $U
\times G$ of $P$, the connection $\nabla$ is given by a $1$-form
$\theta \in \Omega^1(M, \g)$. In the corresponding trivialization
$U \times \g$ of $\adp$, the differential $d_\nabla$ equals $d +
[\theta, \cdot]$. Note that both terms are degree $1$ derivations
of the bracket.
\end{proof}

Recall the definition of the bilinear form $\omega$ from Example
\ref{ex1},
\begin{equation*}\small{
\omega:\Omega^i(M,\adp)\times
\Omega^j(M,\adp)\overset{\wedge}{\to} \Omega^{i+j}(M,\adp\otimes
\adp)\overset{\l\cdot,\cdot\r}{\to} \Omega^{i+j}(M,\mathbb
R)\overset{\int_M}{\to} \mathbb{R}},
\end{equation*}
when $i+j=2d$, and zero otherwise.

It follows from an argument similar to the proof of
\ref{App-compatible-bracket} and the Stokes Theorem
that,

\begin{cor} \label{app-coro-horiztonal}
For $x,y\in\Omega^\ast (M, \adp)$, we have,
$$
\omega(d_{\nabla}x,y)+(-1)^{|x|}\omega(x,d_{\nabla}y)=0.
$$
\end{cor}

\begin{cor} \label{app-coro-d-tr}For a section $s
\in\Gamma(\adu)$, we have, $$d \tr(s)=\tr(\nabla s)\in
\Omega^1(M,\R).$$ Therefore, for $\alpha\in \Omega^\ast (M,\adu)$,
$$d \tr(\alpha)=\tr(d_{\nabla}\alpha).$$
\end{cor}
\begin{proof}
Locally $d_\nabla=d+[\theta, \cdot]$ and $tr [\cdot, \cdot]=0$.
\end{proof}
Let $E, E_1, E_2, ~\text{and}~ E_3$ denote vector bundles over a
compact manifold $M$ of dimension $m$, each endowed with a flat
connection $\nabla$. Given a bundle map $\beta : E_1 \otimes E_2
\to E_3$, there is a natural wedge product $\wedge: \Omega^\ast
(M, E_1) \otimes \Omega^\ast  (M, E_2) \to \Omega^\ast (M, E_3)$.
In addition, if $\beta$ is parallel, this map respects the
differentials and consequently induces a well-defined cup product
$\cup: H^\ast  (M, E_1) \otimes H^\ast  (M, E_2) \to H^\ast (M,
E_3)$ on the cohomology. A noteworthy special case is the
differential graded algebra structure of the differential forms
with values in a flat bundle whose fibres are algebras and the
parallel transport maps are algebra maps. For instance
$\Omega^\ast  (M, \adu)$ is an associative (not graded
commutative) differential graded algebra giving rise to the graded
associative algebra $H^\ast (M, \adu)$. Let us also consider
chains $C_\ast (M, E)$ with values in a bundle $E$. More
precisely, $C_k(M, E)$ is the vector space generated by pairs
$(\sigma, s)$, where $\sigma: \Delta^k \to M$ is a singular chain
and $s$ is a flat section of $\sigma^\ast  E$. The boundary map
$\partial:C_k(M, E) \to C_{k-1}(M, E)$ is defined as $\partial
(\sigma, s)=\Sigma_{i=0}^k (-1)^i( \sigma_i, s|_{\sigma_i})$,
where $\partial \sigma = \Sigma_{i=0}^k (-1)^i \sigma_i$. There is
also a cap product $H^\ast  (M, E_1)\otimes H_\ast  (M, E_2) \to
H_\ast (M, E_3)$.  For example Lemma \ref{lem-vam-vector} uses the
Poincar\'e duality map $H^\ast (M,\adp) \overset{\cap
[M]}{\longrightarrow} H_\ast (M,\adp)$, induced by the bundle map
$\beta : {\bf R} \otimes \adp \to \adp$, defined by $\beta (r,
x)=rx$, where $\R \to {\bf R} \to M$ is the trivial $\R$-bundle.
The Poincar\'e duality map is an isomorphism (see \cite{C}).

\section{Parallel transport and iterated integrals}\label{Appen-Hol}

Let $\mathbb R^k \to E \to M$ be a vector bundle endowed with a
fixed connection $\nabla$ which is not necessarily flat. In this
appendix we show that the holonomy of any other a arbitrary
connection $\nabla'$ (not necessarily flat) along a loop $\gamma$
in $M$ can be expressed by an iterated integral and the holonomy
of $\nabla$. This explains the definition of generalized holonomy
in Section 5. We presume that this is well known but we do not
know a good reference for it.

To start, we trivialize $E$ over a given path $\gamma:[a,b]\to  M$
and using this  trivialization we take $\nabla+\theta\wedge \cdot$
a local expression for $\nabla'$. Here,
$\theta\in\Omega^1(M,End(\mathbb R^k))$ is a matrix valued 1-form.
We aim at finding $t\mapsto P_t \in End(\R^k)$ the parallel
transport along $\gamma$.  For a vector $ v\in \R^k$, $t \mapsto
P_t(v)$ the parallel transport of $v$ along $\gamma$ is defined by
the differential equation,
\begin{equation}\label{eq-parellel}
\nabla'_{\dot{\gamma}(t)} P_t(v)=0.
\end{equation}
 Since the bundle has been trivialized , we can think of $P_t$ as a matrix valued function which satisfies,
 $$
 \nabla_{\dot{\gamma}(t)} P_t(v) +\theta(\dot{\gamma}(t))P_t(v)=0.
$$
 Let $\psi_t$ be the parallel transport of  $\nabla$ and $R_t \in End(\R^k)$ such that $P_t=R_t \psi_t$.
 Since $\nabla_{\dot{\gamma}(t)} \psi_t=0$, by Leibnitz rule the differential equation (\ref{eq-parellel}) reduces
 to,
 \begin{equation}\label{eq-parllel-diff-eq}\frac{dR_t}{dt} +\theta(\dot{\gamma}(t))
 R_t=0,
\end{equation} in the matrix valued functions for the initial value $R_0=id$.
Let $A(t)=\theta(\dot{\gamma}(t))$, the first observation is that
(\ref{eq-parllel-diff-eq}) is equivalent to,
\begin{equation}\label{eq-parllel-diff-eq2}
R_t =id +\int_0^t A(s) R_sds,
\end{equation}
thus $R_t$ is the fixed point of the operator $\phi \mapsto
T(\phi)$, where,
\begin{equation}
 T(\phi)(t)=id +\int_0^t A(s)\phi(s)ds.
\end{equation}
Let $b\in [0,1]$ be such that $R=\int _0^1\parallel A(s)\parallel
ds <1$, then,
\begin{equation}
\parallel T(\phi_1)-T(\phi_2)\parallel\leq R
\parallel\phi_1-\phi_2\parallel,
\end{equation}
hence $T$ is a contracting operator on the Banach space of
continuous matrix valued functions on $ [0,b]$. It is a known fact
that such operators have unique fixed point.  In fact the fixed is
the limit point of the Cauchy sequence $\phi,T(\phi), T^2(\phi),
\cdots $ where $\phi$  is  arbitrary.

By some classical techniques we can show that this  solution can
be extended to the entire interval $ [0,1]$ in order to obtain a
fixed point for the operator $T$  on the Banach space of
continuous matrix valued function on $ [0,1]$.
\begin{prop}
Equation \ref{eq-parellel}  with given initial value has a unique
solution which is differentiable and satisfies
\ref{eq-parllel-diff-eq2}.
\end{prop}
Applying \ref{eq-parllel-diff-eq2} repeatedly we obtain,

\begin{equation}
\begin{split}
R_1 &=id +\int_0^1 A(t_1) R_{t_1}dt_1\\
 &=id +\int_0^1 A(t_1)dt_1+ \int_0^1\int_0^{t_1} A(t_1)A(t_2) R_{t_2}dt_2dt_1\\
 &\phantom{aaaaaaaaaaaaaaaaa}\vdots\\
 &=id +\int_0^1 A(t_1)dt_1+ \int_0^1\int_0^{t_1} A(t_1)A(t_2) dt_2dt_1 \\&+\int_0^1\int_0^{t_1}\int_0^{t_2}
 A(t_1)A(t_2)A(t_3)dt_3dt_2dt_1 + \cdots \\
&+\int_0^1\int_0^{t_1}\cdots \int_0^{t_n} A(t_1)A(t_2)\cdots
A(t_n) dt_n\cdots dt_2dt_1+\cdots.
\end{split}
\end{equation}

Note that since $A(t_i)=i_{\frac{\partial}{\partial t_i}}ev_i^\ast
(\theta)$, we have,
$$
\int_0^1\int_0^{t_1}\cdots \int_0^{t_n} A(t_1)A(t_2)\cdots A(t_n)
dt_n\cdots dt_2dt_1=\int_{\gamma} \int_{\Delta^n} ev_{1}^\ast
(\theta)ev_{2}^\ast (\theta)\cdots ev_{n}^\ast (\theta).
$$
Hence,
\begin{equation}
\begin{split}
R_1 &=\psi_1 +\int_{\gamma} \int_{\Delta^1} ev_{1}^\ast (\theta) \psi_1+ \int_{\gamma} \int_{\Delta^2} ev_{1}^\ast (\theta) ev_{2}^\ast (\theta) \psi_1+ \cdots \\
&+\int_{\gamma} \int_{\Delta^n} ev_{1}^\ast (\theta)ev_{2}^\ast
(\theta)\cdots ev_{n}^\ast (\theta)\psi_1+\cdots.
\end{split}
\end{equation}
\begin{rem} \label{conv} This explains the formula for the generalized holonomy in
section 5 and confirms the convergence of the generalized
holonomy. To sum up, the solution of a general time-dependent
system of linear equations $\dot{X}=A(t)X$ with initial condition
$X(0)=X_0$ may be expressed in terms of a power series each of
whose terms is a Chen iterated integral. Note that in dealing with
a time-independent system, in which $A(t)=A$, for all $t\in \R$,
this formula reduces to $X(t)= e^{tA}X_0$, since the volume of the
$n$-simplex is $\frac{1}{n!}$.
\end{rem}

\end{document}